%%%%%%%%%%%%%%%%%%%%%%%%%%%%%%%%%%%%%%%%%%%%%%%%%%%%%%%%%
%  Lectures on the affine Hecke algebras,		% 
%   and Macdonald conjectures
%   by A. Kirillov, Jr.					%
%%%%%%%%%%%%%%%%%%%%%%%%%%%%%%%%%%%%%%%%%%%%%%%%%%%%%%%%%
% This file should be processed by AmsTeX version 2.1%%%%
%%%%%%%%%%%%%%%%%%%%%%%%%%%%%%%%%%%%%%%%%%%%%%%%%%%%%%%%%
\input amstex
\magnification 1200
\documentstyle{amsppt}
\NoBlackBoxes
\NoRunningHeads
\def\a{\alpha}
\def\l{\lambda}

\def\Z{\Bbb Z}
\def\C{\Bbb C}
\def\R{\Bbb R}
\def\N{\Bbb N}
\def\P{\Cal P}
\def\X{\Cal X}
\def\Y{\Cal Y}
\def\Yt{\widehat {\Cal Y}}
\def\Gt{\widehat {G}}
\def\eps{\varepsilon}
\def\d{\partial}

\def\H{\Cal H}

\def\Hhat{{\hat H}}
\def\Wt{{\widetilde W}}
\def\wt{{\tilde w}}

\def\Rhat{\widehat R}
\def\ahat{{\widehat \alpha}}
\def\lhat{{\widehat\lambda}}
\def\sln{\frak{sl}_n}

\def\Res{\operatorname{Res}}
\def\Ker{\operatorname{Ker}}
\def\v{^\vee}
\def\<{\langle}
\def\>{\rangle}
\def\qbinom#1#2{{\thickfracwithdelims[]\thickness0{#1}{#2}}}

\topmatter
\title Lectures on the affine Hecke algebras and Macdonald conjectures
\endtitle
\author {\rm {\bf Alexander A. Kirillov, Jr.} \linebreak
\vskip .1in
Department of Mathematics\linebreak
Yale University\linebreak 
New Haven, CT 06520, USA\linebreak
e-mail: kirillov\@math.yale.edu\linebreak 
\vskip .1cm 
\bf January 11, 1995}
\endauthor
\thanks The author was supported by Alfred P. Sloan dissertation fellowship 
\endthanks
\endtopmatter

\document

\vskip .5cm
\head Introduction\endhead

These are the notes from a series of lectures the author gave at
Harvard University in the Fall of 1994. The goal of these lectures is
to give a  self-contained exposition of recent
result of Cherednik (\cite{C6}), 
who proved Macdonald's inner product conjectures (see \cite{M2})
for arbitrary root systems, using the algebraic structure he called 
 ``double affine Hecke algebra''.   
I wanted to make these lectures understandable to general
mathematical audience, not only to the experts; for this reason 
  I give all the necessary definitions (including that of
Macdonald's polynomials, Hecke algebras etc) and motivations 
 in the course. Also, I omit some 
technical details of proof (leaving enough information  so that
an experienced reader always can fill the gaps). 
However, I do 
 assume that the reader is familiar with (finite-dimensional and
affine) root systems and Weyl groups.
 A short introduction to these
notions can be found in \cite{Hu1, Chapter III}; see also a recent
survey of Koornwinder \cite{Ko} to see how these notions appear in the
theory of orthogonal polynomials. For more detailed
expositions can we refer the reader to  \cite{B, Hu2, V}.

These lectures are of expository nature; they do not contain any
new results. I include the reference to the original papers in each
chapter. Also, I must stress that these lectures were written under
strong influence of 
(unpublished) lecture course of Ivan Cherednik (Yale, Fall 1991) and
Howard Garland (Yale, Spring 1993) and
series of talks by Ian Macdonald (Yale, October 1993 and Leiden  Univ.,
 May 1994). Large part of my lectures (most of Lectures~2--4,6,7) follows
Macdonald's exposition, so all the credits for these lectures
belong to him; of course,
I am completely responsible for any errors that could be found in
these notes. 

I am deeply grateful to Ivan Cherednik for numerous
discussions, in which he explained to me many parts of this theory.
Without these discussions, my lectures would never come to a happy
end. 

Finally, I would like to thank the mathematics department of Harvard
University for its hospitality during my work on these lectures and
Professors Tom Koornwinder and Masatoshi Noumi for their valuable 
remarks on the preliminary version of these notes.

\newpage
%%%%%%%%%%%%%%%%%%%%%%%%%%%%%%%%%%%%%%%%%%%%%%%%%%%%
\head Lecture 1: Commuting families of differential
operators, Jacobi polynomials and Hecke algebras
\endhead 
%%%%%%%%%%%%%%%%%%%%%%%%%%%%%%%%%%%%%%%%%%%%%%%%%%%%%
This lecture is of introductory nature. In this lecture we
will illustrate the main idea of this course: that for a
certain  family of special functions (Jacobi
polynomials, Macdonald's polynomials), which possess many
interesting special properties (for example, they can be
described as eigenfunctions of a large family of commuting
differential (difference) operators), there exists an algebraic
structure hidden behind them which gives a natural
explanation for these properties. For simplicity, we start
with the classical (differential) case, which is more
geometrical. In fact, complete  proofs
for differential case are more difficult than for the difference one,
but since we are not giving proofs in this
lecture, differential case is quite transparent. 

    We start with the brief survey of the theory of Jacobi
polynomials following the papers of Heckman and Opdam
\cite{HO, H1, O1, O2}. We do not give any proofs; unless
otherwise stated, all the proofs can be found in the above
mentioned papers of Heckman and Opdam (though some of the
results had been known before).
 Let $V$ be a vector space over $\C$,
and $R\subset V$ be a (reduced, irreducible) root system of
rank $n=\text{ dim }V$.
We use the standard notations $(\, ,\,), R^+, \a_i, Q, P, W\ldots$
for the 
inner product in $V$, positive roots, basis of simple roots, root lattice,
weight lattice, Weyl group etc. As usual, we write $\l\le\mu$ if
$\mu-\l\in Q_+$. 
Let us consider the group
algebra of the weight lattice $\C[P]$, which is spanned by
the formal exponentials $e^\l, \l\in P$. We can interpret
them as functions on $V$ by $e^\l(v)=e^{(\l, v)}$. Let us
fix  for every $\a\in R$ a number $k_\a\in \Z_+$ such that
$k_{w(\a)}=k_\a$ for every $w\in W$, and  define the
following differential operator in $V$ (Sutherland operator): 

    $$L_2=\Delta-\sum_{\a\in R^+}k_\a(k_\a-1)
	\frac{(\a,\a)}{(e^{\a/2}-e^{-\a/2})^2},\tag 1.1$$

where $\Delta$ is the Laplace operator. 
Obviously, this operator is $W$-invariant. This operator has many
remarkable properties. For the root system $A_{n-1}$ it appeared as
the Hamiltonian of a system of $n$ particles on a line with potential
of interaction given by $\frac{1}{\sinh ^2 (x_i-x_j)}$ (see \cite{S}).
 This system is
completely integrable (\cite{OP}), which also holds for general root systems: 

\proclaim{Theorem 1.1} {\rm (Heckman, Opdam)}

    Let $\Cal D=\{\text{\rm differential operators in $V$
with coefficients from the ring $\C[P](e^\a-1)^{-1}$}\}$.
Consider $\Bbb D=\{\d\in \Cal D| \d \text{ is
$W$-invariant}, [\d, L_2]=0 \}$. Then $\Bbb D$ is
isomorphic to $(S[V])^W$: there exists ``Harish-Chandra
isomorphism'' $\gamma:\Bbb D\to(S[V])^W$ such that for a
homogeneous $p$, $\gamma^{-1}(p)=p(\d_v)+\text{ lower order
operators}$. \endproclaim

Moreover, for  most root systems (in particular, for $A_n, B_n, D_n$
with $n\ge 4$) it is shown in
 \cite{OOS} that under suitable restrictions  the
only differential operators of the form $\Delta +V(h)$ 
 satisfying this complete integrability
property are Sutherland operator and its rational and elliptic
analogues (with $\frac{1}{\sinh^2 x}$ replaced by $\frac{1}{ x^2}$ and
$\wp (x)$, respectively) and their modifications.

    \demo{Remark} It is relatively easy to construct $\gamma$ and show
that it is injective (and thus,  $\Cal
D\subset (S[V])^W$); the difficult part is to prove that $\gamma$ is
surjective. \enddemo
   
For example, under this isomorphism $\gamma(L_2)=\sum
v_i^2, v_i$ being orthonormal basis in $V$.

    This theorem naturally gives rise to the following
questions: 
\roster
\item What are the eigenfunctions of these operators and
their properties? 

    \item Why is it that $\Bbb D\simeq (S[V])^W$; is there
any natural explanation to this fact as well as
construction of $\gamma^{-1}$?

    \item Is it possible to extend $\gamma^{-1}$ to all
polynomials, i.e. construct for every $v\in V$ a
differential operator $D_v$ such that (1) $D_v$ commute and
(2) $wD_vw^{-1}=D_{w v}$ for every $w\in W$ in such a way
that $\gamma^{-1}(p)=p(D_v)$?
\endroster

    We'll try to answer these questions. Let us start with
the last one. Here the answer is obviously ``no''. It is so
even in the $\frak{sl}_2$ case (i.e., $R$ of type $A_1$),
when $L_2=\frac{d^2}{d x^2}+C\frac{1}{\sinh ^2 x}$, which obviously is
not a square of any first-order differential operator. 

    This analogous to the definition of Dirac operator in
physics. Recall that Dirac operator was introduced  as an
attempt to find a square root of the Laplace operator. Such
a square root does not exist if you look for it in the
class of scalar-valued differential operator. However, if
you consider differential operators with values in the
Clifford algebra then such a square root does exist, and it is
called Dirac operator. 

    Similar construction is possible here, and the
corresponding algebraic structure -- similar to that of
Clifford algebra -- is degenerate affine Hecke algebra.
Before describing it, let us slightly reformulate the
problem. 

Let
$$\gathered
\delta^k=\prod_{\a\in R^+}
(e^{\a/2}-e^{-\a/2})^{k_\a},\\
\rho_k=\frac{1}{2}\sum_{\a\in R^+}k_\a \a.\endgathered\tag 1.2$$

    Now, define

   $$M_2=\delta^{-k}(L_2-(\rho_k, \rho_k))\delta^k.\tag 1.3$$
\proclaim{Proposition 1.2}

    $$M_2=\Delta-\sum_{\a\in R^+}k_\a
\frac{1+e^\a}{1-e^\a}\d_\a.\tag 1.4$$
\endproclaim

    \proclaim{Proposition 1.3} 
\roster \item $M_2$ preserves the space $\C[P]^W$ of Weyl
group invariant polynomials, and so do all the operators
from $\delta^{-k}\Bbb D \delta^k$.  
\item Let $m_\l=\sum_{\mu\in W\l}e^\mu, \l\in P_+$ be the
basis of orbitsums in $\C[P]^W$. Then the action of $M_2$
is triangular in this basis:

    $$M_2 m_\l=(\l, \l+2\rho_k)m_\l+\sum_{\mu< \l}
c_{\l\mu} m_\mu.$$ 
\endroster\endproclaim

    Thus, we can restrict ourselves to considering only the
action of $M_2$ in $\C[P]^W$; instead of considering all
eigenfunctions we consider only symmetric polynomial
eigenfunctions. 

    \proclaim{Definition} Jacobi polynomials $J_\l\in
\C[P]^W, \l\in P_+$ are defined by the following conditions: 
	\roster\item $J_\l=m_\l +\text{lower order terms}.$
    \item $M_2J_\l=(\l, \l+2\rho_k)J_\l.$
\endroster
\endproclaim
It is easy to show that these conditions
define $J_\l$ uniquely. These polynomials have a number of interesting
properties; for example, they are orthogonal with respect to a certain
inner product (we'll discuss it in the next lecture). For special
values of $k$ these polynomials can be interpreted as zonal spherical
functions on a certain Riemannian symmetric spaces associated with the
group $G$ (see \cite{H1});
 in the case of the root system of type $A_n$ they can be
interpreted as zonal spherical functions with the values in a
representation of $G$ for arbitrary $k\in Z_+$ (see \cite{EFK}). 

    Now, let us return to  questions 2 and 3 above: is
it possible to introduce commuting  differential operators
 $D_v, v\in V$ such
that $wD_v w^{-1}=D_{wv}$ for any $w\in W$, and $M_2=\sum
D_{v_i}^2 +const, v_i$ being an orthonormal basis in $V$? As we
have seen before, the answer is ``no''. However, this is
almost possible if we allow $D$ to be not necessarily
scalar differential operators (see explanation of ``almost''
in the remark at the end of this lecture).
 There are two ways in which it can be
done; in fact, they are closely related and can be
considered as special cases of a general approach
(see \cite{C3, C5}), but we won't 
to go into details here. 

    \roster\item We can let $D$ to be not scalar but
matrix-valued differential operators: it is possible to
introduce $D$ acting in some
vector space $E$, and a linear map $E\to \C$ such
that any matrix-valued differential operator which is
obtained as a symmetric polynomial of $D_{v_i}$ can be
pushed forward to some scalar differential operator. In
such a way one can get the commuting family of differential
operators discussed above from the symmetric polynomials in
$D_{v_i}$. This approach was considered in detail by Matsuo
(\cite{Ma}), where he took $E=\C[W]$.  We won't use this
approach in these lectures.

    \item We can let $D$ be scalar valued but not
differential operators: we allow $D$ to include the action
of the Weyl group (which acts on functions by permuting
arguments). These operators are not local; however, if we
take symmetric polynomials of $D_{v_i}$ then these
operators preserve the space of Weyl group invariant
polynomials, and their restrictions on this subspace
coincide with certain differential operators (uniquely defined). In
particular, $\biggl(\sum
D_{v_i}^2\biggr)|_{\C[P]^W}= M_2$. 
\endroster

    We will be mostly interested in this last approach. In
the next lectures we will use it in difference case to get
difference analogue of this commuting family of
differential operators, study their polynomial eigenfunctions
(Macdonald's polynomials) and prove Macdonald's inner 
product identities. Today, we will illustrate the ideas in
a "baby example". Namely, let us consider the rational
degeneration of the above differential operators. Introduce
rescaling operator $A_t$ by $(A_t f)(v)=f(tv)$. Consider 
$L_2(t)=t^{-2}A_t^{-1} L_2 A_t$. Then it is easy to see
that as $t\to 0$, $L_2(t)$ has a limit, which we will call
$L_2^{rat}$: 

    $$L_2^{rat}(v)=\Delta-\sum_{\a\in R^+}
k_\a(k_\a-1)\frac{(\a,\a)}{(\a, v)^2}.\tag 1.5 $$

    Similarly, we can get the rational degeneration of
$M_2$: 

    $$M_2^{rat}(v)=\Delta-\sum_{\a\in R^+}
k_\a\frac{1}{(\a, v)}\d_\a.\tag 1.6$$

    Consider the  $\sln$ case, i.e. the root system of
type $A_{n-1}$. Then $V\subset \C^n$, and we can identify
functions on $V$ with functions on $\C^n$ which are
invariant with respect to the translations $(x_1, \ldots, x_n)\mapsto
(x_1+c, \ldots, x_n+c)$. In this case, $M_2^{rat}$    
takes the form

    $$M_2^{rat}=\sum\d_i^2+2k\sum_{i<j}\frac{1}{x_i-x_j}(\d_i-\d_j),
\tag 1.7 $$
where $\d_i=\frac{\d}{\d x_i}$.

We will show how one can obtain the expression for $M_2^{rat}$ of the
form $M_2^{rat}=\sum D_i^2$ using so-called differential-difference
operators introduced by Dunkl (see \cite{H2}). 
Let  $b_{ij}$ be the following operator:

    $$b_{ij}f=\frac{s_{ij}f-f}{x_i-x_j},\tag 1.8$$
where $s_{ij}$ acts on functions of $x_1, \ldots, x_n$ by
permutation of  arguments: $x_i\leftrightarrow x_j$. Note
that $b_{ij}$ preserves the space $\C[x_1, \ldots, x_n]$,
since $x^ny^m-x^m y^n$ is divisible by $x-y$. Also, it is easy to see
that $wb_{ij}w^{-1}=b_{w(i) w(j)}$ for any $w\in S_n$. 

    Define (rational) Dunkl operators by
$$D_i=\d_i-k\sum_{j\ne i}b_{ij}\tag 1.9$$
(sometimes they are also called local Dunkl operators, as opposed to
operators with trigonometric coefficients, which are called ``global'').

\proclaim{Theorem 1.4} 
\roster\item $wD_i w^{-1}=D_{w(i)}$ for any $w\in S_n$. 
\item $[D_i, D_j]=0$.
\endroster
\endproclaim
\demo{Proof} (1) is trivial. Let us prove (2). Obviously,
$[D_i, D_j]$ contains only terms linear in $k$ and
quadratic in $k$. Since $[\d_i, b_{lm}]=0$ if $i\ne l,m$,
the term linear in $k$ is equal to $-k([\d_i, b_{ji}]-[\d_j,
b_{ij}])$. Since $b_{ij}=-b_{ji}$, this equals
$k[\d_i+\d_j, b_{ij}]=0$. 

    Standard arguments show that in order for the term
quadratic in $k$ vanish it is necessary and sufficient that
$b_{ij}$ satisfy the classical Yang-Baxter equation: 

    $$[b_{12},b_{13}]+[b_{12},b_{23}]+[b_{13},b_{23}]=0.$$
This can be proved by direct calculation, which is rather
boring. We will show another way to prove it later. 
\qed\enddemo

Now, for every operator of the form $D=\sum_{w\in S_n} D_w w,
D_w$ -- a differential operator with rational coefficients, define the
associated differential operator by 

$$\Res\bigl(\sum D_w w\bigr)=\sum D_w.\tag 1.10$$

Note that if $D$ preserves the space of symmetric polynomials then 
so does $\Res D$, and $D|_{\C[P]^W}=\Res D|_{\C[P]^W}$. 

\proclaim{Theorem 1.5} $\sum D_i^2$ is $S_n$-invariant, preserves 
$\C[x_1, \ldots, x_n]^{S_n}$ and $\Res \biggl(\sum
D_i^2\biggr)=M_2^{rat}$. \endproclaim

\demo{Proof} Invariance of $(\sum D_i^2)$ and the
fact that it preserves $\C[x_i]^{S_n}$ follow from the previous theorem. To
calculate $\Res (\sum D_i^2)$, write 

    $$\sum D_i^2=\sum_i\left(\d_i^2 
-k\sum_{j\ne i}(\d_ib_{ij}+b_{ij}\d_i)
+4k^2\sum\Sb j\ne i\\ l\ne i\endSb b_{ij}b_{il}\right).$$

    Since $b_{ij}|_{\C[x]^{S_n}}=0$ (this is the crucial
step!), we have $\Res \sum D_i^2=\Res
(\sum_i(\d_i^2-k\sum_{j\ne i} b_{ij}\d_i))$. It is easy to
check that 
$$b_{ij}\d_i=(x_i-x_j)^{-1}(s_{ij}-1)\d_i=(x_i-x_j)^{-1}(\d_j
s_{ij} -\d_i)$$
and thus

    $$\Res \bigl(\sum D_i^2\bigr)=\sum _i \left(\d_i^2-k\sum_{j\ne i}
(x_i-x_j)^{-1}(\d_j-\d_i)\right) =  M_2^{rat}.$$
\qed\enddemo

\proclaim{Corollary} $M_2^{rat}$ can be included in 
 a commutative family of symmetric differential operators: 
 $M_1^{rat}=\sum
\d_i, M_2^{rat}, \ldots,  M_n^{rat}$ with coefficients from $\C[x_1,
\ldots , x_n](x_i-x_j)^{-1}$. 
\endproclaim

\demo{Proof} Take $M_r^{rat}=\Res \sum D_i^r$. \enddemo

    Thus, we have proved (in this baby example) the complete
integrability  theorem~1.1 
 and gave an explicit construction of these
differential operators. However, there are some questions even in this
case, namely:
\roster\item Where did we get these expressions for $D_i$ and $b_{ij}$
from? Is there a way to guess them?

\item Why do $b_{ij}$ satisfy Yang-Baxter equation? 
\endroster

The answer to these questions is that there exists some simple
algebraic construction which allows to get both the expression for
$b_{ij}$ and their properties without any calculations. This is the
degenerate affine Hecke algebra. 

\proclaim{Definition} The degenerate affine Hecke algebra  for the root
system $A_{n-1}$ is the algebra $H'_n$ over $\C$ spanned by its two
subalgebras $\C[S_n]$ and $\C[x_1, \ldots, x_n]$ with the relations 

$$\gathered x_{i+1}s_i-s_ix_i=s_i x_{i+1} - x_i s_i=h\\
  x_is_j=s_j x_i\quad\text{ if } i\ne j, j+1,\endgathered\tag 1.11$$
where $h\in \C$ is a fixed constant, and $s_i=s_{i, i+1}, i=1\ldots
n-1$ are the standard generators of $S_n$. 

\endproclaim 

This algebra is a deformation of the semidirect product
$\C[S_n]\ltimes \C[x]$, which is its limit as $h\to 0$
 (for brevity, we write $\C[x]$ for $\C[x_1, \ldots, x_n]$).

\proclaim{Theorem 1.6} $H'_n=\C[S_n]\cdot \C[x]=\C[x]\cdot \C[S_n]$,
i.e. every $g\in H'_n$ can be uniquely written in every of the
following forms: 

$$g=\sum_{w\in S_n}p_w(x) w=\sum_{w\in S_n} w q_w(x),\tag 1.12$$ 
where $p_w, q_w\in \C[x]$. \endproclaim

We do not prove this theorem here, referring the reader to \cite{C5}
and references therein. Note that it is obvious that every element can
be written in either of the forms in (1.12); the difficult part is to
prove the uniqueness. 

\remark{Remark} Of course, one could as well define degenerate affine
Hecke algebra for any root system, and the analogue of Theorem 1.6 is
also true; see the above cited paper of Cherednik. \endremark

Now, let $E$ be any module over $S_n$. Define the induced module $\hat
E=\operatorname{Ind}_{\C[S_n]}^{H'_n} E$. It follows from the theorem
above that as a linear space (and, moreover, as a $\C[x]$-module),
$\hat E\simeq E\otimes \C[x]$. In particular, let us take $E=\C$ with
the trivial action of $S_n$. Then
$\hat E\simeq \C[x]$, and we get the following proposition:

\proclaim{Proposition 1.7} There exists a unique action of $H'_n$ in
$\C[x]$ such that 
\roster\item $p\in \C[x]\subset H'_n$ acts by multiplication by $p$.
\item $\hat s_i 1=1$, where $\hat s$ denotes the action of an element $s\in
S_n\subset H'_n$ in $\C[x]$. 
\endroster
\endproclaim

\proclaim{Proposition 1.8} In the above defined representation

$$\hat s_i=s_i+ h b_{i, i+1},\tag 1.13$$
where $s_i$ is the usual action of $S_n$ in $\C[x]$ \rom{(}by permutation
of $x_i$\rom{)}, and $b_{ij}$ is defined by \rom{(1.8)}.\endproclaim

\demo{Proof} We know from Proposition 1.7 that $\hat s_i$ are defined
uniquely. Moreover, it is easy to see that in fact they are defined
uniquely by the commutation relations (1.11) and the condition $\hat
s_i 1=1$, so there is no need to check that $\hat s_i$ satisfy the
relations of the symmetric group. But  $\hat s_i$ defined by
formula (1.13) satisfy both $\hat s_i 1=1$ (obvious) and (1.11), which
can be shown by a rather short explicit calculation.\qed\enddemo

Thus, we see that the operators $b_{ij}$ we defined before have a very
natural interpretation in terms of the degenerate affine Hecke algebra: they
describe the action of it in the induced representation. Now, let us
show that this allows to prove the classical Yang-Baxter equation
without any calculations. Indeed, it follows from Propositions 1.7 and
1.8 that $\hat s_i$ defined by (1.13) satisfy the braid relations: 

$$\hat s_{1}\hat s_{2} \hat s_{1}=\hat s_{2} \hat s_{1}\hat
s_{2}.$$ 

Let us define $R_{12}=s_1 \hat s_1, R_{23}=s_2 \hat s_2$, and
$R_{13}=s_1 R_{23} s_1 =s_2 R_{12}s_2$ (check this last identity!).
It is easy to check  that the braid relation for
$\hat s_i$ implies  (quantum) Yang-Baxter equation for $R_{ij}$:

$$R_{12} R_{13}R_{23} = R_{23} R_{13} R_{12}.$$
 Since $R_{ij}=1+hb_{ij}$, it is a standard fact that
(quantum) Yang-Baxter equation for $R_{ij}$ implies classical
Yang-Baxter equation for $b_{ij}$. 

  Thus, we have shown that using the notion of the degenerate affine
Hecke algebra along with the ``Poincare-Birkhoff-Witt theorem'' 1.6,
we can easily and naturally derive the formula for $D_i$ and prove all
the required properties. This is the main idea of this course. In
the next lectures we will explain in detail how similar technique works
in the difference case. 

\remark{Remark} In fact, this technique requires a bit more
careful approach even in the trigonometric differential case. Namely,
in this case we can not construct $D_v$ which would commute and
satisfy the relation $w D_v w^{-1}=D_{wv}$.  However, we can construct
something very close to it: we can construct $D_v$ which commute and
their commutation relations with $S_n$ are given by  formula (1.11),
i.e. they satisfy the relations of degenerate affine Hecke algebra.
This still allows us to get a commuting family of differential
operators, since  it is known that $\C[x]^{S_n}$ is the center of $H'_n$.
\endremark

\newpage
%%%%%%%%%%%%%%%%%%%%%%%%%%%%%%%%%%%%%%%%%%%%%%%%%%%%%%%%%
\head Lecture 2: Macdonald's polynomials and difference
operators \endhead
%%%%%%%%%%%%%%%%%%%%%%%%%%%%%%%%%%%%%%%%%%%%%%%%%%%%%%%%%%%
	
	Now we start a systematic study of the difference case. We
won't use the first lecture (except as motivation).

    In this lecture we define the quantum analogue of the
Jacobi polynomials discussed last time. Unless otherwise stated, 
the results in this
lecture are due to Macdonald (\cite{M1, M2}). 

    We begin with fixing the notations. Let $V$ be a
finite-dimensional vector space over $\R$ with a
positive definite symmetric bilinear form $(\cdot, \cdot)$,
$R\subset V$ be a reduced irreducible 
root system. We fix a decomposition of
$R$ into positive and negative roots: $R=R_+\sqcup R_-$ and
denote by $\a_1, \ldots , \a_n$ the basis of simple roots
in $R$. For every root $\a$ define the dual root
$\a\v=\frac{2\a}{(\a, \a)}$. Denote by $Q=\bigoplus \Z \a_i$
the root lattice, $Q_+=\bigoplus \Z_+\a_i$, $P=\{\l\in V |
(\l, \a_i\v)\in \Z\}$ the weight lattice, and $P_+=\{\l\in
V| (\l, \a_i\v)\in \Z_+\}$ the set of dominant integral
weights. It has a natural basis of fundamental weights
$\omega_i: (\omega_i, \a_j\v)=\delta_{ij}$. In a similar
way, define $Q\v=\bigoplus \Z\a_i\v$ be the coroot lattice,
$Q_+=\bigoplus \Z_+\a_i\v$, $P\v=\{\l\in V|(\l, \a_i)\in
\Z\}$ -- coweight lattice, $P_+\v=\{\l\in V|(\l, \a_i)\in
\Z_+\}$ -- dominant coweights, $b_i$ -- fundamental
coweights: $(b_i, \a_j)=\delta_{ij}$. As usual, we define 
the highest root $\theta\in R$ by $\theta-\a\in Q_+$ for
all $\a\in R$, and the element $\rho=\frac{1}{2}\sum_{\a\in
R_+}\a$; then $(\rho, \a_i\v)=1$, so $\rho\in P$. 

    For every $\a\in R$, let $s_\a$ be the corresponding
reflection: $s_\a(v)=v-(v, \a_i\v)\a_i$, and let $W$ be the
Weyl group generated by $s_\a$. Our main goal is
construction of a special basis in the space $\C[P]^W$. 
One example of such a basis is given by the orbitsums: 
$m_\l=\sum_{\mu\in W\l}e^\mu$. This basis is orthogonal with
respect to the following inner product in $\C[P]$:
$\<f,g\>_0= \frac{1}{|W|}[f\bar g]_0$, where the bar involution is
defined by $\overline{e^\l}=e^{-\l}$, and $[\,\,]_0$ is the
constant term: $[\sum a_\l e^\l]_0=a_0$.

    We will generalize this basis as follows. 
    Suppose that for every $\a\in R$ we have a variable
$t_\a$ subject to the conditions $t_\a=t_{w(\a)}$ (so, we
have at most 2 different $t$). Consider the field of
rational functions in $t_\a$ and one more independent
variable $q$: $\C_{q,t}=\C(q, t_\a)$. Define the inner
product in $\C_{q,t}[P]$ by 
	$$\<f, g\>_{q,t}=\frac{1}{|W|}[f\bar g\Delta_{q,t}]_0\tag 2.1$$ 
where bar involution is extended by $\C_{q,t}$ linearity,
and 

    $$
\Delta_{q,t}=\prod_{\a\in R}\prod_{i=0}^\infty
	\frac{1-q^{2i}e^\a}{1-t_\a^2q^{2i}e^\a}
	=\Delta^+_{q,t}\overline{\Delta^+_{q,t}},
    \tag 2.2$$
where 
$$\Delta^+_{q,t}=\prod_{\a\in R^+}\prod_{i=0}^\infty
	\frac{1-q^{2i}e^\a}{1-t_\a^2q^{2i}e^\a}\tag 2.3$$

Both $\Delta$ and $\Delta^+$ should be considered as Laurent series in
$q, t$ with coefficients from $\C[P]$; then the inner product also
takes values in Laurent series. It is easy to see that  it is non-degenerate
 and $W$-invariant. If $t_\a=1$ for all $\a$ then this inner
product coincides with previously defined $\<\cdot,
\cdot\>_0$. 

    \proclaim{Theorem 2.1}{\rm (Macdonald)} There exists a
unique family of elements $P_\l\in \C_{q,t}[P]^W$, $\l\in P_+$
satisfying the following conditions:
\roster\item $P_\l=m_\l+\sum_{\mu<\l}a_{\l\mu}m_\mu$

    \item $\<P_\l, P_\mu\>_{q,t}=0$ if $\l\ne \mu$
\endroster\endproclaim

These polynomials are called Macdonald's polynomials (see
\cite{M1, M2}). They form a basis in $\C_{q,t}[P]^W$. Note
that the theorem above is not trivial: since $<$ is not a
complete order, you can not get $P_\l$ by orthogonalization
of $m_\l$. 

\demo{Examples} \roster\item 
If $t_\a=1$ then independently of $q$, $P_\l=m_\l$.

    \item If $t_\a=q$ for all $\a$ then $P_\l=\chi_\l$ are  Weyl
characters.

	\item If $q, t\to 1$ in such a way that $t_\a=q^{k_\a},
k_\a\in \Z_+$ are fixed then $P_\l\to J_\l$, where $J_\l$ are Jacobi
polynomials defined in the previous lecture as eigenfunctions of some
differential operator $M_2$. Indeed, in this limit $\Delta^+_{q,t}\to
\delta^k$ and it suffices to check that $M_2$ is self-adjoint with
respect to the inner product $\<f, g\>=\frac{1}{|W|} [f\bar
g(\delta\bar\delta)^k]_0=\<f\delta^k, g\delta^k\>_0$. This is equivalent
to the fact that $L_2$ is self-adjoint with respect to $\<\, ,\,\>_0$,
which is obvious.  

    \item If $q=0$, $t_\a=1/p, p\in \Z_+$ is a prime then
$P_\l$ can be interpreted as zonal spherical functions on
the corresponding $p$-adic group. 
	
	\item For the root system $A_1$ Macdonald's polynomials
coincide with so-called $q$-ultraspherical polynomials (\cite{AI}),
 which are a special case of Askey-Wilson polynomials (\cite{AW}).
More generally, Macdonald's polynomials for (non-reduced) root system
$BC_1$ are precisely Askey-Wilson polynomials.

\endroster\enddemo

    For simplicity, we only consider in these lectures the case when 
$t_\a=q^{k_\a}$ for $k_\a\in \Z_+$. It is not an important
restriction: all the results that we prove  can
 be generalized to independent $q,t$ without much problem. However,
this case is easier from the technical point of view; for example, 
formula (2.2) takes the form 

    $$\Delta_{q,t}=\prod_{\a\in R}\prod_{i=0}^{k_\a-1}
(1-q^{2i}e^\a)$$
thus allowing to avoid infinite products. We
will use the notations $\<\, ,\,\>_k, \Delta_k$ for $\<\, ,\,\>_{q,t},
\Delta_{q,t}$ etc. Also, it will be convenient to use 
$\rho_k=\frac{1}{2}\sum_{\a\in R^+} k_\a \a$.
We always consider $q$ as a formal variable; the results 
also hold if $q$ is a complex number provided that it is
not a root of unity.

Similarly to the classical case, for some special values of $k$
Macdonald's polynomials can be interpreted as zonal spherical functions
on certain $q$-symmetric spaces associated with the group $G$ (see
\cite{N}).

    \demo{ Sketch of proof of theorem 2.1} To prove the
theorem, it suffices to 
find an operator $D:\C_q[P]^W\to
\C_q[P]^W$ such that 
\roster \item
$D m_\l=\sum\limits\Sb\mu\in P_+\\ \mu\le \l\endSb c_{\l\mu}m_\mu$

    \item $c_{\l\l}$ are distinct. 
    
    \item $D$ is self-adjoint with respect to the inner product
$\<\, , \>_{q.t}$.

    \endroster

Let us construct such an operator. 
\proclaim{Definition} A weight $\l\in P^+$ is called
minuscule if $0\le (\l, \a\v)\le 1$ for every $\a\in
R^+$.\endproclaim 

    It is easy to see that a weight $\l\ne 0$ can be minuscule
only if $\l$ is one of the fundamental weights. Indeed, let
$\phi$ be the highest root for the root system $R\v$;
then $\phi=\sum n_i \a_i\v, n_i\in \N$. Thus, $(\l,
\phi)\le 1$ implies that $\l=0$ or $\l=\omega_r$ for some $r$ such
that $n_r=1$. In fact, the following is known: 

    \proclaim{Lemma 2.2}{\rm (see \cite{B,V})} The set of
all minuscule weights 
is a system of representatives for $P/Q$, i.e. every $\l\in P$ can be
written in a unique way 
in the form $\l=b+\a$ for some minuscule weight $b$ and $\a\in Q$. \endproclaim

    This implies that there are no non-zero minuscule weights for
the root systems $E_8, F_4, G_2$ and that the number of
minuscule weights is always not greater than the rank of
the root system with the equality only for $A_n$. 

    Let $\pi\in P\v$ be a minuscule coweight: $0\le (\pi,
\a)\le 1$ for $\a\in R^+$. Define
$T_\pi:\C_q[P]\to \C_q[P]$ by $T_\pi(e^\l)=q^{2(\l,
\pi)}e^\l$ (this requires adding appropriate
fractional powers of $q$ to $\C_q$). Let $D_\pi$ be defined
by 

    $$D_\pi(f)=\sum_{w\in W}
w\biggl(\frac{T_\pi(\Delta^+f)}{\Delta^+}\biggr)
=\sum_{w\in W} 
	w\biggl(\prod\Sb \a\in R^+\\ (\a,\pi)=1\endSb 
		\frac{1-q^{2k_\a}e^\a}{1-e^\a}T_\pi(f)
	\biggr) 
   \tag 2.4$$

    Let us prove that $D_\pi f\in \C_q[P]^W$ for every
$f\in \C_q[P]$. It is obvious from (2.4) that $D_\pi f$ 
is $W$-invariant rational function with poles only on $e^\a
-1=0$ and all the poles are simple. Thus, $\delta D_\pi f$
(where $\delta=\prod_{\a\in R^+} (e^{\a/2}-e^{-\a/2})$ is the Weyl
denominator) is a $W$-antiinvariant element of $\C_q[P]$. As is
well-known, this implies that in fact $D_\pi f\in
\C_q[P]^W$. 

    The triangularity condition (1) above can be easily
verified by direct calculation, which also shows that 
$c_{\l\l}=q^{2(\pi, \rho_k)}\sum_{w\in W} q^{2(\pi,
w(\l+\rho_k))}$ (which requires the identity
$\rho_k-w(\rho_k)=\sum\limits_{\a\in R^+\cap w^{-1}R^-} k_\a \a$).

    The self-adjointness of $D_\pi$ can be easily deduced  from the
definition of the inner product. Finally, one can check that
for all root systems having non-zero minuscule coweights, except $D_n$,
one can find a minuscule coweight $\pi$ such that the
corresponding eigenvalues $c_{\l\l}$ are distinct; for
$D_n$ it is not so, but there exists a linear
combination of operators $D_\pi$ corresponding to 
 minuscule coweights such that the eigenvalues
are distinct. This proves the theorem for all cases when
non-zero minuscule coweights exist, i.e. all cases except $E_8, F_4,
G_2$.\enddemo

    The above proof used that $\pi$ is a minuscule weight:
otherwise, you could get a product of factors of the form
$(1-e^\a)(1-q^2e^\a)\ldots $ in the denominator. Thus, this
proof fails if there are no non-zero minuscule coweights (i.e., for
the root systems $E_8, F_4, G_2$). This can be fixed,
which, however, requires certain ingenuity; we refer the
reader to the original papers of Macdonald.

    \proclaim{Proposition 2.3} If $\pi_1, \pi_2$ are
minuscule coweights then $D_{\pi_1}, D_{\pi_2}$
commute.\endproclaim 

    Proof of this proposition is straightforward. 

    \demo{Example} Let $R\subset V$ be of type $A_{n-1}$;
we identify $V$ with a subspace in $\R^n$ given by the
condition $\sum \l_i=0$. Then $e^\l\mapsto x_1^{\l_1}\ldots
x_n^{\l_n}$, where $\l=(\l_1, \ldots, \l_n)\in V\subset
\R^n$ gives an isomorphism of $\C_q[P]$ with the space of
homogeneous polynomials in $x_i^{\pm 1}$ of degree zero.
In this case all the fundamental weights (=coweights)
$\omega_1, \ldots, \omega_{n-1}$ are minuscule, and the
corresponding operators take the form

    $$D_r=\sum_w w\bigl(
		\prod\Sb i,j\\ i\le r< j \endSb
	\frac{x_i-q^{2k}x_j}{x_i-x_j} T_{1\ldots r}\bigr)
	=r!(n-r)!\sum\Sb I\subset\{1 \ldots n\}\\ |I|=r\endSb
      \left(	
	\prod\Sb i\in I \\ j\notin I\endSb 
	\frac{x_i-q^{2k}x_j}{x_i-x_j}
      \right)
	T_I,\tag 2.5$$
where $T_I=\prod_{i\in I}T_i, (T_i f)(x_1, \ldots,
x_n)=f(x_1, \ldots, q^2x_i, \ldots, x_n)$. 

    \enddemo

    In this example, we have constructed a commuting family
of difference operators, and the number of independent 
operators is equal to the rank of the root system, which is
the natural quantum analogue of the commuting families of
differential operators considered in Lecture~1. However,
this analogy fails for other root systems, since the number
of difference operators we get from the minuscule coweights
is in general less than the rank of the root system. The
correct answer is that there always exists a commuting
family of difference operators, and their number is equal
to the rank of the root system, but they are not obtained
from minuscule coweights. One of the main goals of the next
lectures will be construction of these difference
operators based on the representation theory of affine Hecke
algebras. In general, the explicit expressions for these
operators are rather complicated, which makes them impossible
to guess; the fact that the operators corresponding to
the simple coweights can be written by such a simple
formula as (2.4) is a lucky exception. Explicit expressions can also
be written for the
non-reduced root system $BC_n$, though they are more complicated (see
\cite{D}). 

    Now we can formulate the Macdonald's inner product
identity. 

    \proclaim{Theorem 2.4} {\rm (Macdonald's inner product
identity)} 
$$\aligned
\<P_\l, P_\l\>_{q,t}=&\prod_{\a\in R^+}
\prod_{i=1}^{k_\a-1} 
\frac{1-q^{2(\a\v, \l+\rho_k)+2i}}
    {1-q^{2(\a\v, \l+\rho_k)-2i}}\\
=& q^{\sum_{\a\in R^+}k_\a (k_\a-1)}
\prod_{\a\in R^+}
\prod_{i=1}^{k_\a-1} 
\frac{[(\a\v, \l+\rho_k)+i]}
    {[(\a\v, \l+\rho_k)-i]},\endaligned\tag 2.6$$
where $[n]=\frac{q^n-q^{-n}}{q-q^{-1}}$. \endproclaim

    This formula has been conjectured by Macdonald, who has
proved it for the root systems $A_n$ (unpublished). Also, it has been
proved for the (not reduced) root system  $BC_1$ (see \cite{AW}).
 The first  proof for arbitrary root system 
was given by Cherednik; we will give this proof in the
following lectures. 

    \demo{Examples}
\roster\item Let $k_\a=1$. Then $P_\l=\chi_\l$ are Weyl
characters, and (2.6) reduces to $\<\chi_\l, \chi_\l\>_1=1$,
which is well-known. 

    \item Let $\l=0$. Then $P_\l=1$, and (2.6) reduces to 
$$\frac{1}{|W|}[\Delta_{q,t}]_0
	= q^{\sum_{\a\in R^+}k_\a (k_\a-1)}
	\prod_{\a\in R^+}\prod_{i=1}^{k_\a-1}
\frac{[(\a\v, \rho_k)+i]}{[(\a\v, \rho_k)-i]},\tag 2.7$$
which is known as constant term identities.
These identities were first conjectured by Macdonald
(\cite{M3}), though some special cases had been  known before.
Classical ($q=1$) case  of these identities was proved by Opdam
(\cite{O3}), using the technique of shift operators, which we will
discuss later.
In general case, these 
identities had been proved case-by-case 
for most root systems, with the exception
of $E$ series (see \cite{BZ}, \cite{GG}, \cite{Ha}, \cite{K}).

    These identities  can be
rewritten in the following form:
$$\left[\prod_{\a\in R^+}\prod_{i=0}^{k_\a-1}
(1-q^{2i}e^\a)(1-q^{2i+2} e^{-\a})\right]_0 
=\prod\qbinom{kd_i}{k},\tag $2.7'$ $$
where $\qbinom{a}{b}=\frac{[a]!}{[b]![a-b]!}$ is the
$q$-binomial coefficient, and $d_i$
are exponents of the Weyl group $W$, i.e. the
degrees of the free generators of $(S[V])^W$. This
reformulation is not  trivial: it involves some identity for the
Poincar\'e series of $W$, which can be found in 
\cite{M4}. In particular, for $A_{n-1}$ the
exponents are $2,3,\ldots, n$, and  formula ($2.7'$) becomes
$$\left[\prod _{i<j}\prod_{l=0}^{k-1}
	\biggl(1-q^{2l}\frac{x_i}{x_j}\biggr)
	\biggl(1-q^{2l+2}\frac{x_j}{x_i}\biggr)
   \right]_0=
    \frac{[nk]!}{[k]!^n}.\tag 2.8$$

    It is also worth noting that if we let $k\to \infty$ in
($2.7'$) then the formula we get is closely related with 
the denominator identity for the
corresponding affine root system, also due to Macdonald. 
\endroster
\enddemo

\newpage
%%%%%%%%%%%%%%%%%%%%%%%%%%%%%%%%%%%%%%%%%%%%%%%%%%%
\head Lecture 3: Affine Hecke algebras and induced
representations\endhead 
%%%%%%%%%%%%%%%%%%%%%%%%%%%%%%%%%%%%%%%%%%%%%%%%%%%

    In this lecture we define and study the affine Hecke
algebras; later the results of this lecture will be used to
obtain Macdonald's difference operators. Most of the results we give
in this lecture are due to Lusztig (\cite{L}); in less general case
they were first proved by Bernstein and Zelevinsky (unpublished). 

	We start with the definition of affine Weyl group. 
In the notations of previous sections, let $\hat V=V\oplus
\C\delta$; we will interpret elements of $\hat V$ as
functions on $V$ by $(v+k\delta)(v')=(v, v')+k$.  Define the
affine root system $\Rhat=R\times \Z\delta$ and the
positive affine roots by $\Rhat^+=\{\a+k\delta| \a\in R,
k>0 \text{ or } \a\in R^+, k\ge 0\}$. The basis of simple
roots in $\Rhat^+$ is given by $\a_0=-\theta+\delta, \a_1,
\ldots, \a_n$. Equivalently, $\Rhat^+=\{\ahat\in
\Rhat|\ahat \text{ is non-negative on }C \}$, where $C=\{v\in
V|  \a_i(v)\ge 0, i=0\ldots n\}$ is the (affine) Weyl
chamber. 

For every $\ahat=\a+k\delta$ we define the reflection 
$s_\ahat:\hat V\to \hat V$ by 

    $$s_\ahat:\lhat\mapsto \lhat -(\l, \a\v)\ahat,\tag
3.1$$
where $\lhat=\l+m\delta, \ahat=\a+k\delta$. Note that this
action preserves $\Rhat$. 

    Dual action of $s_\ahat$ in $V$ is just the reflection
in $V$ with respect to the (affine) hyperplane
$\ahat(v)=0$: 

    $$s_\ahat: v\mapsto v - \ahat(v) \a\v\tag 3.2$$
     
    We will will use the notations $s_0\ldots s_n$ for
$s_{\a_0}\ldots s_{\a_n}$. Define the affine Weyl group
$W^a$ as the group  generated by
$s_\ahat$. Then the following facts are well-known:

    \proclaim{Theorem 3.1}
\roster\item $W^a=W\ltimes \tau(Q\v)$, where the action of
$\a\v\in Q\v$ in $\hat V$  is given by 

    $$\tau(\a\v):\lhat\mapsto \lhat - (\a\v, \l)\delta \tag
3.3$$ 
\rom{(}Respectively, the dual action in $V$ is given by
$\tau(\a\v):v\mapsto v+\a\v$.\rom{)}

In particular, $s_0=\tau(\theta\v)s_\theta=s_\theta\tau(-\theta\v)$. 

    \item $W^a$ is generated by $s_0,\ldots, s_n$ with the
relations $s_i^2=1$ and the  Coxeter relations: for $i\ne j$  

    $$s_i s_j s_i \ldots = s_j s_is_j \ldots \quad\text{
$m_{ij}$ terms on each side}, \tag 3.4 $$
where $m_{ij}$ are defined in the standard way from the
corresponding affine Dynkin diagram. \rom{(} For the root
system $\hat A_1$ there are no Coxeter relations, which is
sometimes formally expressed by letting $m_{ij}=\infty$.\rom{)}

    \item For every $w\in W^a$, its length $l(w)$ with respect
to the generators $s_0,\ldots s_n$ is equal to 

    $$l(w)=|\Rhat^+\cap w^{-1}\Rhat^-	|\tag 3.5$$

    \item $C$ is a fundamental domain for action of $W^a$
in $V$. \endroster\endproclaim

    Now, let us define extended Weyl group $\Wt$ as the
semidirect product $\Wt=W\ltimes \tau(P\v)$, where the
action of $\tau(P\v)$ is given by the same formulas as we
had before for $\tau(Q\v)$. {\bf Note that the action of
extended Weyl group in $\hat V$ preserves $\Rhat$.} It is
easy to see that $W^a$ is a normal subgroup in $\Wt$ and 
$\Wt/W^a\simeq P\v/Q\v$ is an abelian group, whose elements
are in one-to-one correspondence with minuscule coweights
(see Lemma~2.2 in the previous lecture). 
It turns out that $\Wt$ can be presented as a semidirect
product. Let us define the length $l(\wt)$ for arbitrary
$\wt\in \Wt$ by the same formula (3.5). In general, $\Wt$
is not a Coxeter group, and $l(\wt)$ can not be interpreted
as a length of a reduced decomposition. Moreover, you can
have elements of length~0. 
Define $\Omega =\{\wt\in \Wt|l(\wt)=0\}=\{\wt\in
\Wt| \wt(C)=C\}$. Obviously, this is a subgroup, and it
follows from the theorem above that $\Wt=\Omega\ltimes
W^a$; thus, $\Omega \simeq P\v/Q\v$. This means that every
element of $\Omega$ has the form $\pi_r=\tau(b_r)w_r$ for
some minuscule coweight $b_r$ and $w_r\in W^a$. It is also
useful to note that $\pi_r$ acts on the simple roots
$\a_0\ldots \a_n$ by some permutation; in particular,
$\pi_r(\a_0)=\a_r$. Thus, we get the following description
of $\Wt$: 

    $$\gathered
\Wt=\Omega\ltimes W^a,\\
\text{with the relation }\pi_r s_i\pi_r^{-1}=s_j\quad 
\text{if } \pi_r(\a_i)=\a_j.\endgathered\tag 3.6$$

    We will need some properties of the length function $l(\wt)$.
\proclaim{Lemma 3.2} 
\roster

\item $l(\pi_r \wt)=l(\wt)$

    \item 
$l(\wt s_i)= \cases l(\wt)+1, &\wt(\a_i)\in \Rhat^+\\ 
l(\wt)-1, &\wt(\a_i)\in \Rhat^-\endcases$ 

\item If $w\in W, \l\in P\v$ then 
$$l(w\tau (\l))=\sum_{\a\in R^+} |(\l, \a)+\chi(w\a)|,\tag
3.7$$ 
where $\chi(\a)= 0$ if $\a\in R^+$ and $1$ if $\a \in R^-$.
\endroster\endproclaim

\demo{Proof} (1) is obvious from the definition; (2) follows
from (1) and standard results about affine Weyl group; (3) can
be derived straightforwardly from the definition of length.
\enddemo

    \proclaim{Corollary 3.3}
\roster\item If $\l\in P\v$ then $l(\tau(\l))=2(\l^+, \rho)$, where
$\l^+$ is the dominant coweight lying in the $W$-orbit of $\l$. 

\item If $\l\in P\v_+$ then
$l(w\tau(\l))=l(w)+l(\tau(\l))$. 

\item If $(\l, \a_i)=0, i=1,\ldots, n$ then
$l(\tau(\l)s_i)=l(s_i\tau(\l))=l(\tau(\l))+1$. 

\item If $(\l, \a_i)=-1$ then
$l(s_i\tau(\l))=l(\tau(\l))-1$.  

    \endroster\endproclaim

    Now we can define the braid group.
\proclaim{Definition 3.4} The braid group $B$ is the group
generated by the elements $T_\wt, \wt\in \Wt$ modulo the
following relations: 
$$T_v T_w=T_{vw}\quad \text{ if } l(vw)=l(v)+l(w), \quad v,w \in\Wt.
\tag 3.8$$
\endproclaim

    In particular, this implies that the elements $T_u,
u\in \Omega$ form a subgroup in $B$ which is isomorphic to
$\Omega$; abusing the language, we will use the same
notation $\pi_r$ for $T_{\pi_r}$. Also, we will write $T_i$
for $T_{s_i}, i=0\ldots n$. 

\proclaim{Lemma 3.5} 
$T_{\wt s_i}=T_{\wt} T_i$ if $\wt s_i$ is a reduced expression, and
$T_{\wt s_i}=T_\wt T_i^{-1}$ otherwise.\endproclaim
\demo {Proof} Immediately follows from Lemma 3.2(2).\enddemo

    \proclaim{Theorem 3.6} $$B=\Omega\ltimes B(T_0, T_1,
\ldots, T_n),$$
where $B(T_0, \ldots, T_n)$ is the group with generators
$T_0, \ldots, T_n$ and relations \rom{(3.4)} \rom{(}Coxeter
relations\rom{)},
and the action of $\Omega$ on $T_i$ is given by 
$\pi_rT_i \pi_r^{-1}=T_j$ if $\pi_r(\a_i)=\a_j$. 
\endproclaim
\demo{Proof}This theorem immediately follows  from
the previous results and the following well-known result, 
due to Iwahori and Matsumoto:
for any two reduced expressions for an element $\wt\in W^a$
one can be obtained from
anoher by a sequence of Coxeter relations, i.e. without
using the relations $s_i^2=1$. \enddemo

    Now we can make one of the most crucial steps. Define
elements $Y^\l\in B$ for $\l\in P\v$ by
\roster\item If $\l\in P\v_+$ then  $Y^\l=T_{\tau(\l)}$.

    \item If $\l=\mu-\nu, \mu, \nu\in P\v_+$ then
$Y^\l=Y^\mu(Y^\nu)^{-1}$. 
\endroster

\proclaim{Theorem 3.7} \roster\item
$Y^\l$ is well-defined for all $\l$, and
$Y^\l Y^\mu=Y^{\l+\mu}$. 
\item If $\tau(\l)=\pi_r s_{i_l}\ldots s_{i_1}$ 
is a reduced expression, denote 
$\a^{(1)}=\a_{i_1}, \a^{(2)}=s_{i_1}(\a_{i_2}), \ldots,
\a^{(l)}=s_{i_1}\ldots s_{i_{l-1}}(\a_{i_l})$ the associated sequence
of affine roots \rom{(}we'll discuss them in detail in
Lecture~\rom{5)}.  Then
 $Y^\l=\pi_r T_{i_l}^{\eps_l}\ldots T_{i_1}^{\eps_1}$, where 
$\eps_i=1$ if the corresponding $\a^{(i)}$ has the form 
$\a^{(i)}=\a+k\delta, \a\in R^+$, and $\eps_i=-1$ otherwise. 
 In particular, if $\l\in P\v_+$ then
all $\eps_i=1$ and if $\l \in P\v_-$ then all $\eps_i=-1$. \endroster
\endproclaim
\demo{Proof} (1) is trivial in view of
$T_{\tau(\l+\mu)}=T_{\tau(\l)} T_{\tau(\mu)}$ if $\l,
\mu\in P_+\v$.

As for (2), the proof consists of several steps which we briefly
outline. Let us call an affine root $\ahat$ $R$-positive if 
$\ahat=\a+k\delta, \a\in R^+$. 
  \proclaim{Lemma} For any \rom{(}not necessarily reduced\rom{)}
 expression $\tau(\l)=\pi_r s_{i_l}\dots
s_{i_1}$ denote $\tilde Y^\l=\pi_r T_{i_l}^{\eps_l}\ldots
T_{i_1}^{\eps_1}$, where the signs $\eps_i$ are determined as in the
statement of the theorem. Then $\tilde Y^\l$ does not depend on the choice
of the expression for $\tau(\l)$.\endproclaim

    \demo{Proof} First, show that if $\tau(\l)=x s_i^2 y =
x y, x, y\in \Wt$ are two expressions for $\tau(\l)$ then
the corresponding expressions for $\tilde Y^\l$ are equal.
Indeed, the sequence of roots associated with the first
expression differs from the second one by insertion of the
pair $\a^{(k)}=y^{-1} \a_i, \a^{(k+1)}=y^{-1}(-\a_i)$.
Since precisely one of these roots is $R$-positive,
we have $T_i^{\eps_{k+1}}T_i^{\eps_k}=1$. 

Similarly, if $\tau(\l)=x (s_i s_j)^{m_{ij}} y =
x y, x, y\in \Wt$, where $m_{ij}$ is as in (3.4), 
 are two expressions for $\tau(\l)$ then
the sequence of roots associated with the first
expression differs from the second one by insertion of
the set of roots $\{y^{-1}\beta_k\}_{\beta_k\in R<\a_i,
\a_j>}$,  where $R<\a_i,
\a_j>$ is the root system of rank two spanned by $\a_i,
\a_j$. Moreover, these roots appear in their natural cyclic
order: one can choose an orientation in $\R \a_i\oplus
\R\a_j$ such that the roots $\beta_k$  appear in the counterclockwise
order. Both of these facts can be easily checked
case-by-case, since one only has to consider root systems
of rank 2. Since the condition of being an $R$-positive root 
specifies a half-space, we see that among the roots
$y^{-1}\beta_k$ precisely one half (i.e. $m_{ij}$) is 
$R$-positive, and they go in a row. Thus, the corresponding
part of the expression for $\tilde Y^\l$ has the form 

$$\undersetbrace p\text{ ``$+$'' signs} \to{T_iT_j T_i\dots}\quad
 \undersetbrace m_{ij}\text{ ``$-$'' signs} 
		\to {\dots T_i^{-1}T_j^{-1}\dots}
     \undersetbrace m_{ij}-p\text{ ``$+$'' signs} \to{\dots T_iT_j }$$

    or

$$\undersetbrace p \text{ ``$-$'' signs}
	\to{T_i^{-1} T_j^{-1} T_i^{-1}\dots}\quad
\undersetbrace m_{ij}\text{ ``$+$'' signs} 
	\to {\dots T_i T_j \dots}
\undersetbrace m_{ij}-p \text{ ``$-$'' signs}
	\to{\dots T_i^{-1}T_j^{-1} }$$
for some $p$. In both cases, this product is equal to 1,
which completes the proof of the Lemma.\enddemo

    Now it is relatively easy to prove (2). Indeed, let
$\l\in P_+\v$, and let $\tau(\l)=\pi_r s_{i_l}\ldots
s_{i_1}$ be a reduced expression. Then it is known (see,
for example, \cite{Hu2}) that 
$\{\a^{(1)}, \dots, \a^{(l)}\}=\{\ahat\in
\Rhat^+|\tau(\l)\ahat\in \Rhat^-\}$. Explicit calculation
shows that all of them are $R$-positive, and therefore
$\tilde Y^\l=\pi_r T_{i_l}\dots T_{i_1}=Y^\l$. Finally, it
follows from the fact that $\ahat$ is $R$-positive iff
$\tau(\l)\ahat$ is $R$-positive and the Lemma proved above that 
$\tilde Y^\l \tilde Y^\mu = \tilde Y^{\l+\mu}$. 
Thus, $\tilde Y^\l=Y^\l$ for all $\l\in P\v$, which
 concludes proof of (2).
\qed\enddemo

Thus, the subgroup generated by $Y^\l$ is
isomorphic to coweight lattice $P\v$. 
    
\proclaim{Lemma 3.8} 
\roster\item The elements $Y^\l, \l\in P\v, T_1, \ldots,
T_n$ generate $B$ as a group. 
\item If $(\l, \a_i)=0, 1\le i\le n$ then $T_iY^\l =Y^\l
T_i$. 

    \item If $(\l, a_i)=1, 1\le i \le n$ then $Y^\l=T_i
Y^{s_i\l} T_i$. \rom{(}Note: this is {\bf not} a misprint: there are two $T_i$
and no $T_i^{-1}$ in the formula.{)}

\endroster\endproclaim

    \demo{Proof} 
\roster \item Since every $\wt\in \Wt$ can be written as
$\wt=\tau(b_r)w$ for some minuscule (and thus, dominant) 
weight $b_r$ and $w\in
W^a$, it follows from Lemma 3.5 that $Y^\l, T_0, T_1,\ldots
T_n$ generate $B$. Similarly, it follows from
$s_0=\tau(\theta\v)s_\theta$ that $T_0$ can be written in
terms of $Y^{\theta\v}$ and $T_1, \ldots, T_n$. 

    \item If $\l\in P\v_+$ this follows immediately from
Corollary 3.3. Result for general $\l$ follows from the
fact that every $\l$ can be presented in the form
$\l=\mu-\nu, \mu, \nu\in P\v_+, (\mu, \a_i)=(\nu, \a_i)=0$.

    \item Suppose $\l\in P\v_+$. Introduce $\pi=\l+s_i(\l)=
2\l-\a_i\v$; then $\pi\in P\v_+$ (obvious). If
$l(\tau(\l))=2(\l, \rho)=p$ then $l(\tau(\pi))=2p-2$, and
it follows from the definition that $l(\tau(\l)s_i)=p-1$. Thus,
if we write $s_i\tau(\pi)=(\tau(\l)s_i) (\tau(\l))$ (it is
easy to see that it is in fact an identity in $\Wt$), then
both left-hand and right-hand sides are reduced
expressions, and thus 
$T_i Y^\pi=T_{\tau(\l)s_i}Y^\l=Y^\l T_i^{-1}Y^\l$, which is
equivalent to the desired equality. 

    If $\l\notin P\v_+$, it can be written as $\l=\mu-\nu,
\mu,\nu\in P\v_+, (\mu, \a_i)=1, (\nu, \a_i)=0$, and thus
the statement follows from previous arguments. 
    
\endroster\enddemo

    Now we are ready for the main definition of this
lecture. Suppose that for every $\a\in \Rhat$ we have a
variable $t_\a$ such $t_\a=t_{w(\a)}$ for every $w\in \Wt$
(thus, you can have at most two different variables). 
Let $\C_t=\C(t_\a)$ be the field of rational functions
in $t_\a$. 

    \definition{Definition} Affine Hecke algebra $\Hhat$ is
the quotient of the group algebra $\C_t[B]$ by the ideal
generated by the following relations: 
	$$(T_i-t_i)(T_i+t_i^{-1})=0, i=0,\ldots, n,\tag 3.9$$
where $t_i=t_{\a_i}$ \rom{(}in particular,
$t_0=t_{\a_0}=t_{\theta}$\rom{)}. 
\enddefinition

    Note that these relations imply that
$T_i^{-1}=T_i+(t_i^{-1}-t_i)$. 

    In a similar way, we can define $H^a$ as a subalgebra
of $\Hhat$ generated by $T_0, \ldots, T_n$ and non-affine
Hecke algebra $H$ as a subalgebra generated by $T_1,
\ldots, T_n$ , so $H\subset H^a\subset \Hhat$ (which is a
complete analogue of $W\subset W^a\subset \Wt$). 

    \proclaim{Theorem 3.9} $\Hhat=\Omega\ltimes H^a$, where
the action of $\pi_r$ on $T_i$ is the same as in Theorem
\rom{3.6}.\endproclaim 

    \proclaim{Lemma 3.10} One has the following relations in
$\Hhat$: 

    $$T_iY^\l - Y^{s_i(\l)}T_i= (t_i-t_i^{-1})
\frac{s_i -1}{Y^{-\a_i\v} -1} Y^\l,\quad i=1, \ldots, n\tag 3.10$$
\rom{(}from now on, expressions of the form $\frac{s_i
-1}{Y^{-\a_i\v}-1}$ stand for $(Y^{-\a_i\v}-1)^{-1} (s_i
-1)$. It is easy to see that the right-hand side is in fact
a polynomial, i.e. lies in $\C_t[Y]$\rom{)}.
    \endproclaim

\demo{Proof} First, simple calculation shows  that if this 
relation is true for $Y^\l, Y^\mu$ then it is also true for
$Y^{\l+\mu}$. Thus, it suffice to prove (3.10) when $(\l,
\a_i)=0$ or $(\l, \a_i)=1$. If $(\l, \a_i)=0$, (3.10)
reduces to $T_iY^\l =Y^\l T_i$, which is Lemma~3.8(2).
Similarly, if $(\l, a_i)=1$ then
$Y^{s_i(\l)}=Y^{\l-\a_i\v}$, and (3.10) reduces to 
$T_iY^\l-Y^{s_i(\l)}T_i=(t_i-t_i^{-1})Y^\l$, which is an
immediate corollary of Lemma 3.8 and identity 
$T_i=T_i^{-1}+(t_i-t_i^{-1})$.\enddemo

\proclaim{Theorem 3.11} \roster
\item
$\Hhat=H\cdot\C_t[Y]$, where the commutation relations of
$H$ and $Y$ are given by \rom{(3.10)}.

    \item Every element of $\Hhat$ can be uniquely
written in any of the following forms: 

$$h=\sum_{w\in W} p_w(Y)T_w=\sum_{w\in W} T_wq_w(Y)\tag 3.11$$
\endroster\endproclaim

    \demo{Proof} One direction is rather easy. 
First, it follows from Lemma 3.8 that the
elements $T_w, w\in W$ and $Y^\l, \l\in P\v$ generate
$\Hhat$ as an algebra. Now relations (3.10) imply that
every element $h$ can be written in the form $h=\sum _wT_w
q_w(Y)$, which proves the existence part of the theorem.

To prove uniqueness, we must show   that these elements are
independent. This is much more difficult, and we do not give the proof
here, referring the reader to \cite{Hu2, Chapter 7}.
\enddemo

    \proclaim{Theorem 3.12} The center $\Cal
Z(\Hhat)=\C_t[Y]^W$ \rom{(}i.e., the elements from $\C_t[Y]$
which are invariant with respect to the 
usual action of $W$\rom{)}.\endproclaim

    \demo{Proof} It is easy to check, using (3.10), that
$\C_t[Y]^W\subset \Cal Z(\Hhat)$. On the other hand, in the
specialization $t_i=1$ we have $\Cal Z(\Hhat)=\C[Y]^W$,
which is easy to prove; thus, the same must be true for
generic $t_i$. \enddemo

    Now, we will discuss the representations of $\Hhat$.
General theory of representations of $\Hhat$ is rather
complicated (see \cite{L}); however, we will use only
representations of some special form. Namely, let $E$ be an
arbitrary representation of $H$; it is known that for
general values of $t_i$ $H$ is isomorphic as an algebra to
the group algebra $\C_t[W]$, and thus has the same
representations. Define a representation $\hat E$ of
$\Hhat$ as an induced representation: 
$\hat E=\text{Ind}_H^\Hhat E$. It follows from Theorem~3.11
that as a vector space, $\hat E=\C_t[Y]\otimes E$,
and the action of $\C_t[Y]$ is by left multiplication. In
particular, let us take the trivial representation of $H$,
i.e. let $E=\C_t, T_i\mapsto t_i, i=1, \ldots, n$. Then we
get an action of $\Hhat$ in the space $\C_t[Y]$. 

    \proclaim{Theorem 3.13} The above defined action of
$\Hhat$ in $\C_t[Y]$ is given by 
$$\gathered Y^\l\mapsto Y^\l\\
T_i\mapsto t_is_i+(t_i-t_i^{-1})\frac{s_i-1}{Y^{-\a\v_i}-1}.\endgathered
\tag 3.12$$
\endproclaim
\demo{Proof} Immediately follows from Lemma 3.10.
\enddemo

    \remark{Remark} The degenerate affine Hecke algebra
which we defined in Lecture~1 (for the root system
$A_{n-1}$) can be obtained as the following degeneration of 
the affine Hecke algebra $\Hhat$: write (formally)
$Y^\l=t^{y_\l/h}$. Then as $t\to 1$ the relations for $y_\l,
T_i$ become the relations of the degenerate affine Hecke
algebra. 
   \endremark 

\newpage
%%%%%%%%%%%%%%%%%%%%%%%%%%%%%%%%%%%%%%%%%%%%%%%%%%%%%%%%%%%%%%%
\head Lecture 4: Double affine Hecke algebras and commuting difference
operators\endhead 
%%%%%%%%%%%%%%%%%%%%%%%%%%%%%%%%%%%%%%%%%%%%%%%%%%%%%%%%%%%%

So far, we have just explained the preliminaries. Now we are ready to
prove some of the results which we promised in Lecture~2. Unless
otherwise stated, all constructions and results in this lecture are
due to Cherednik.

Let us recall some facts from the last lecture. We have defined affine
Hecke algebra $\Hhat$ which has two descriptions: 

\roster\item $\Hhat$ is generated by elements $\pi_r\in \Omega, T_0,
\ldots, T_n$ with relations 
\itemitem{(a)} Coxeter relations for $T_i, i=0,\ldots, n$
\itemitem{(b)} $(T_i-t_i)(T_i+t_i^{-1})=0, i=0,\ldots, n$
\itemitem{(c)}  $\pi_rT_i\pi_r^{-1}=T_j$ if $\pi_r(\a_i)=\a_j$

\item $\Hhat=H\cdot\C_t[Y]$, where $H$ is the Hecke algebra generated by
$T_1, \ldots, T_n$ with the relations above, and $\C_t[Y]$ is the
algebra generated by $Y^\l, \l\in P\v$ such that $Y^{\l+\mu}=Y^\l
Y^\mu$ (thus, $\C_t[Y]\simeq \C_t[P\v]$), and the commutation
relations between $T_i, Y^\l$ are 

$$\aligned 
	T_i Y^\l=Y^\l T_i\quad & \text{ if } (\l, \a_i)=0\\
	T_i Y^\l-Y^{s_i(\l)}T_i= (t_i-t_i^{-1})Y^\l\quad &\text{ if } 
			(\l, \a_i)=1\endaligned$$
\endroster

Also, we have proved that 

$$ \gathered
Y^\l\mapsto Y^\l\\
T_i\mapsto t_is_i+(t_i-t_i^{-1})
\frac{s_i-1}{Y^{-\a_i\v}-1}\endgathered$$

gives a representation of $\Hhat$ in $\C_t[Y]$. 

Now, we want to use this algebra (along with the above defined
representation) to construct a commuting family of difference
operators, in analogy with what we have done in Lecture 1 for
classical case. Note, however, that in Lecture~1 we had to use the
operators $\d_i=\frac{\d}{\d x_i}$ {\bf and} the degenerate affine Hecke
algebra to construct commuting differential operators; we can say that
we have used the algebra spanned by $\d_i, x_i$ and $\hat s_i$. In
this case, it was not necessary since the commutation relations of
$\d_i$ with $\hat s_i, x_i$ were quite simple. However, it turns out
that in order to construct the difference operators, we must use the
quantum analogue of this latter algebra, and not only the part
generated by $x_i, \hat s_i$. 

Such an analogue was constructed by Cherednik, who called 
 ``the double affine Hecke
algebra''. This is an algebra which is generated by three sets of
variables: 

\roster\item $T_i, i=1\ldots n$
\item $Y^\l, \l \in P\v$
\item $X^\mu, \mu\in P$\endroster

Here $Y^\l, T_i$ must satisfy the relations of affine Hecke algebra
above. Since $W$ is the Weyl group for $R\v$ as well as for $R$, we
can define the relations between $X^\mu$ and $T_i$ to be the relations of
the affine Hecke algebra for the root system $R\v$, i.e. the same
relations as above with $Y^\l$ replaced by $X^\mu, \mu\in P$ and
$\a_i$ replaced by $\a_i\v$. What is difficult to define are the
relations between $Y^\l, X^\mu$. For this reason, let us use the
representation (1) above for the affine Hecke algebra $\Hhat^Y$ 
generated by $Y^\l, T_i, i=1,\ldots, n$. Thus, instead of describing
relations between $Y^\l$ and $X^\mu$ we have to describe relation
between $T_0, \pi_r$ and $X^\mu$. 

    Let us consider the affine weight lattice $\hat
P=\{\mu +k\delta |\mu \in P, k\in \frac 1m \Z\}\subset \hat
V$, where  $m\in \Z_+$ is such that
$(\l,\mu)\in \Z$ for every $\l\in P\v, \mu\in P$. Then we
have a natural action of the extended Weyl group
$\Wt=W\ltimes \tau (P\v)$ in $\hat P$ and thus action
of it in the group algebra $\C[\hat P]$ spanned by
$X^{\hat\mu}$. Let us denote $X^\delta= q^{-2}$, where $q$
is an independent variable. Then we have inclusion:
$\C[\hat P]\subset\C_q[P]$, where $\C_q$ is the field of
rational functions in $q^{1/m}$. The action of $\Wt$ can be
extended by $\C_q$-linearity to $\C_q[P]$; in particular, 
$$\gathered
\tau(\l) X^\mu = q^{2(\l,\mu)} X^\mu\\
s_0 X^\mu = X^\mu (X^\theta q^2)^{-(\mu, \theta\v)}.
\endgathered
$$

    Now we are ready for the main definition in this
lecture (and probably in the whole course):

    \proclaim{Definition 4.1} Double affine Hecke algebra $\H$
is an algebra over the field $\C_{q,t}$ of rational
functions of $q^{1/m}, t_\a$ which is generated by elements
$\pi_r\in \Omega, T_0,\ldots, T_n, X^\mu$, $ \mu\in P$ subject
to the following relations:
\roster\item The relations \rom{(a)--(c)} of the affine Hecke
algebra between $T_i, \pi_r$.

    \item $X^{\mu +\nu}=X^\mu X^\nu$

    \item $$\gathered T_iX^\mu=X^\mu T_i\quad \text{ if }
(\mu, \a_i\v)=0\\
T_i X^\mu-X^{s_i(\mu)}T_i=(t_i-t_i^{-1})X^\mu\quad
\text{ if } (\mu, \a_i\v)=1\endgathered\tag 4.1$$

    Here $i=0\ldots n$, where by definition
$\a_0\v=-\theta\v$.

\item $\pi_r X^\mu\pi_r^{-1}=X^{\pi_r(\mu)}$. 
\endroster\endproclaim

    Note that relation (4.1) for $i=0$ reads 
$T_0 X^\mu- X^{\mu +\theta}q^2 T_0=(t_0-t_0^{-1})X^\mu$ if
$(\mu, \theta\v)=-1$.

    Obviously, the subalgebra generated by $\pi_r,
T_0\ldots T_n$ satisfies the relations of  the affine Hecke
algebra  defined in the previous lecture (and in fact, is
isomorphic to it, though we haven't proved this so far). Thus, we can
define the elements $Y^\l, \l\in P\v$ in $\H$. Similarly,
the subalgebra generated by $T_1, \ldots, T_n, X^\mu$
satisfies the relations of affine Hecke algebra for the
root system $R\v$. We will denote these affine Hecke
subalgebras by $\Hhat^Y$ and $\Hhat^X$ respectively. 
In fact, one can check that in the above definition $X^\mu$
and $Y^\l$ play symmetric role: $\H$ could be as well
defined as an algebra spanned by $Y^\l, T'_0, T_1, T_n,
\pi_r'\in \Omega\v$ (see details in \cite{C7}). It is also
worth noting that it is rather difficult to write down
explicitly the commutation relations between $X^\mu$ and
$Y^\l$. For the root system of the type $A_{n-1}$ it can be
described in the topological language: in this case $\H$ is
a deformation of  the braid group of $n$ points on a torus  
factored by the additional relations
$(T_i-t_i)(T_i+t_i^{-1})=0$. Under this correspondence,
$X_i$ corresponds to $i$-th point going around the
$x$-cycle on the torus, $Y_i$ corresponds to $i$-th point
going around the $y$-cycle on the torus, and $T_i$
corresponds to the transposition of the $i$-th and $i+1$-th points
(see \cite{C1, Definition~4.1}).

    \proclaim{Theorem 4.2} Every element $h\in \H$ can be
uniquely written in the form 
$$\sum_{\l, \mu, w} a_{\l\mu w}X^\mu Y^\l T_w,\qquad  a_{\l\mu
w}\in \C_{q,t}\tag 4.2$$
\endproclaim

    \demo{Proof} The existence of such a representation is
quite a standard exercise; the uniqueness is highly
non-trivial, and we postpone the proof until next lecture 
(see Corollary~5.8).\enddemo

\proclaim{Theorem 4.3} The following formulas give a representation of
$\H$ in $\C_{q,t}[X]$: 

$$\gathered 
	\pi_r\mapsto \pi_r\\
	T_i\mapsto t_is_i+(t_i-t_i^{-1})\frac{s_i-1}{X^{-\a_i}-1},
\quad i=0, \ldots, n
\endgathered $$
\endproclaim
\demo{Proof} It turns out that there is nothing to prove: all the
identities we have to check involve at most two $T_i$, and since every
pair of vertices in the affine Dynkin diagram is belongs to some
 subdiagram of finite type (with the exception of the root system $\hat
A_1$),  this theorem follows from the similar
statement for affine Hecke algebra, which we discussed in the previous
lecture. The case of $\hat A_1$ can be easily checked by direct
calculation.  \enddemo

\demo{Example} In this representation, $T_0$ acts as follows:

$$T_0: X^\mu\mapsto \left(t_0 (X^\theta q^2)^{-(\mu, \theta\v)}
+(t_0-t_0^{-1})\frac{(X^\theta q^2)^{-(\mu, \theta\v)}-1}{X^\theta q^2
-1}  \right)X^\mu.\tag 4.3$$

\enddemo

It turns out that in fact this representation is faithful (we'll prove
it later); we will identify elements of $\H$ with the corresponding
operators in $\C_{q,t}[X]$. 

It is clear that for every $\wt\in \Wt$ the action of the
corresponding operator $T_\wt$ can be written as 

$$T_\wt=\sum\Sb \l\in P\v \\w\in W\endSb g_{\l, w} \tau(\l)w\tag 4.4$$
for some $g_{\l, w}\in \C_{q,t}[X](X^\a-1)^{-1}$ 
(recall that $\tau(\l)$ acts in
$\C_{q,t}[X]$ by $X^\mu\mapsto q^{2(\l, \mu)}X^\mu$). In particular,
the same is true for $Y^\l\in \H$. Let $f\in \C [Y]^W=\Cal
Z(\Hhat^Y)$. 

\proclaim{Lemma 4.4} The operators $f(Y)$ preserve the space
$\C_{q,t}[X]^W$. \endproclaim

\demo{Proof} The proof is based on the following simple observation:
$p\in \C_{q,t}[X]$ is $W$-invariant if and only if $T_i p= t_i p$ for
all $i=1,\ldots, n$, which immediately follows from the formula  for
the action of $T_i$. Now, let $p\in \C_{q,t}[X]$. Then
$T_i f(Y)p =f(Y)T_i p=t_i f(Y) p$ (since $\C[Y]^W$ is the center of
$\Hhat^Y$), and thus $f(Y)p\in \C_{q,t}[X]$.\enddemo

Now, for every operator of the form (4.4) define its restriction by 

$$\Res\left(\sum g_{\l, w} \tau(\l)w\right)=\sum g_{\l, w} \tau(\l)
.\tag 4.5$$

This definition is chosen so that (1) $\Res D$ is a difference operator
(that is, it only involves rational functions of $X$ and operators
$\tau(\l)$, not the action of $W$) and (2) if $D$ preserves
$\C_{q,t}[X]^W$ then so does $\Res D$, and $(\Res D)|_{\C_{q,t}[X]^W}
=D|_{\C_{q,t}[X]^W}$.

Then Lemma 4.4 immediately implies:

\proclaim{Theorem 4.5} The operators $L_f=\Res f, f\in
\C_{q,t}[Y]$ commute and are $W$-invariant. \endproclaim

\demo{Proof} First, $L_f|_{\C[X]^W}=f|_{\C[X]^W}$. Thus, $L_f L_g p=
L_g L_f p$ for every symmetric polynomial $p$. Since both $L_f, L_g$
are difference operators (they do not contain the action of the Weyl
group), it is a well-known fact that this implies  $L_fL_g
=L_gL_f$. \enddemo

Thus, we have constructed a commutative family of $W$-invariant 
difference operators in $\C_{q,t}[X]$, labeled by $f\in
\C_{q,t}[Y]^W$. The main goal of the following lectures
will be to show the relation of these operators with the
theory of Macdonald's polynomials. In particular, we will
show that this family includes the Macdonald's difference
operators,  constructed in Lecture~2 for
minuscule weights, and that the eigenfunctions of these
operators are Macdonald's polynomials. 

\demo{Example} Let $R$ be of type $A_1$. In this case there is only
one positive root $\a$, and only one minuscule (co)weight $\rho=\a/2$. The
reduced expression for $\tau(\rho)$ is $\tau(\rho)=\pi_\rho s_1$,
where $\pi_\rho = \tau(\rho)s_1$ is the element of  zero length; it
acts on simple affine roots by permuting $\a_1=\a$ and
$\a_0=-\a+\delta$. In this case, $Y^\rho=T_{\tau(\rho)}=\pi_\rho T_1,
Y^{-\rho} = T_1^{-1}\pi_\rho = \pi_\rho T_0^{-1}= \pi_\rho (T_0+
(t^{-1}-t))$. 

In this case $\C_{q,t}[X]$ is just the  polynomials in $X^{\pm 1/2}$,
where $X=X^\a$, and the action of extended Weyl group is  given by 

$$\gathered\tau (\rho) X=q^2 X\\
\pi_\rho X= q^{-2}X^{-1},\endgathered$$
so the action of the corresponding affine Hecke algebra is given by

$$\gathered
Y^\rho=\pi_\rho\left(ts_1+(t-t^{-1})\frac{s_1-1}{X^{-1}-1}\right)
 =\tau(\rho)\left(t+(t-t^{-1})\frac{1-s_1}{X-1}\right)\\
Y^{-\rho} = \pi_\rho\left(ts_0+(t-t^{-1})\frac{s_0-1}{q^2X-1}
			   +(t^{-1}-t)\right) \\
=\tau(\rho)s_1\left(t\tau(\a)s_1
		+(t-t^{-1})\frac{\tau(\a)s_1-q^2X}{q^2X-1} \right)\\
=t\tau(-\rho)
+(t-t^{-1})\tau(\rho)\frac{\tau(-\a)-q^2X^{-1}s_1}{q^2X^{-1}-1}\\
=t\tau(-\rho) +(t-t^{-1})\frac{\tau(-\rho)-X^{-1}\tau(\rho) s_1} 
			      {X^{-1}-1}	
\endgathered$$

Thus, 

$$\gathered
\Res Y^\rho =t\tau(\rho)\\
\Res Y^{-\rho}=\frac{t X^{-1}-t^{-1}}{X^{-1}-1} \tau(-\rho)
	+(t-t^{-1})\frac{1}{X-1}\tau(\rho),\endgathered
$$

so

$$\Res (Y^\rho+Y^{-\rho})=
\frac{tX-t^{-1}}{X-1}\tau(\rho) +\frac{t
X^{-1}-t^{-1}}{X^{-1}-1}\tau(-\rho),$$
    which is nothing but Macdonald's difference operator
$D_1$ for the root system $A_1$ (cf. formula (2.5)), multiplied by
$t^{-1}$. 
\enddemo

\newpage
%%%%%%%%%%%%%%%%%%%%%%%%%%%%%%%%%%%%%%%%%%%%%%%%%%%%%%%%%    
    \head Lecture 5: Macdonald's difference operators
from double affine Hecke algebras. \endhead
%%%%%%%%%%%%%%%%%%%%%%%%%%%%%%%%%%%%%%%%%%%%%%%%%%%%%%%%%

    Let us recall some facts from the last lecture. We have
defined the double affine Hecke algebra $ \H$ which is
generated by the elements $T_w, w\in W, Y^\l, \l\in P\v,
X^\mu, \mu\in P$. Also, we have defined its representation
in the space $\C_{q,t}[X]$, where $X^\mu$ acts by
multiplication, and $T_i, i=0\ldots n$ act by 
$$T_i\mapsto t_is_i+(t_i-t_i^{-1})\frac{s_i-1}{X^{-\a_i}-1}.$$

    In this lecture we will establish the connection
between this construction and Macdonald's theory; the
variables $ q, t_\a$ used in the definition of the double
affine Hecke algebra will be identified with Macdonald's
parameters $q,t_\a$.

Let us rewrite the expression for $T_i$ as follows: 

$$T_i= s_i G(\a_i),$$
where

$$\aligned
G(\a)&=t_\a+(t_\a-t_\a^{-1})\frac{1-s_\a}{X^\a-1}\\
&=\frac{t_\a X^\a-t_\a^{-1}}{X^\a-1} 
	-(t_\a-t_\a^{-1})\frac{s_\a}{X^\a-1}.
\endaligned\tag 5.1$$

Thus defined $G(\a)$ satisfy 

$$\wt G(\a)\wt^{-1}=G(\wt(\a)).$$

Using this, we can rewrite  for arbitrary $\wt\in \Wt$
  the action of  $T_\wt$ as follows

    $$T_\wt=\wt G(\a^{(l)})\ldots G(\a^{(1)}),$$
where 
 $\a^{(i)}$ are defined from a reduced expression for
$\wt$: if $\wt=\pi_rs_{i_l}\ldots s_{i_1}$ is reduced then
let $\a^{(1)}=\a_{i_1}, \a^{(2)}=s_{i_1}(\a_{i_2}), \ldots,
\a^{(l)}=s_{i_1}\ldots s_{i_{l-1}}(\a_{i_l})$. It is known
that for so defined $\a^{(i)}$, we have 

    $$\{\a^{(1)}, \ldots \a^{(l)}\}=R_\wt:= \Rhat^+\cap
\wt^{-1} \Rhat^-.$$

    We will also need the expressions for $Y^\l$. Recall
(see Theorem~3.7(2)) that if $\tau(\l)=\pi_r s_{i_l}\ldots
s_{i_1}$ is a reduced expression then $Y^\l=\pi_r
T_{i_l}^{\eps_l} \ldots T_{i_1}^{\eps_1}$ for some choice
of signs $\eps_i\in \{\pm 1\}$. Since
$T_i^{-1}=T_i+(t_i^{-1}-t_i)$, this implies that 

$$Y^\l=\tau(\l)G^\pm(\a^{(l)})\ldots G^\pm(\a^{(1)})$$
for some choice of the signs $\pm$, where $G^+(\a)=G(\a)$
and 

    $$G^-(\a)=G(\a)-(t_\a-t_\a^{-1})s_\a
	=\frac{t_\a X^\a-t_\a^{-1}}{X^\a-1} 
	-(t_\a-t_\a^{-1})\frac{X^\a s_\a}{X^\a-1}.\tag 5.2$$

    These expressions for $T_\wt, Y^\l$ are rather complicated
because the expression for $G^\pm (\a)$ is a sum of two terms,
one of which contains $s_\a$. The main goal of today's
lecture is to define some notion of "leading term" of
$T_\wt$ in such a way that the terms with $s_\a$ from
(5.1), (5.2) (or at least as many of them as possible)
would not contribute to the leading term thus making it
easy to compute. 

We start with the definition of a new order on $P\v$. 

    \proclaim{Definition 5.1} Let $\l, \mu\in P\v$. We write
$\l\prec \mu$ if
\roster 
\item $\l^+<\mu^+$, where $\l^+$ is the dominant coweight
lying in the orbit of $\l$, and similarly for $\mu^+$, 

    or 

   \item $\l^+=\mu^+$ and $\l>\mu$ \rom{(}note the change of
sign!\rom{)}.
\endroster 
\endproclaim

    Note that it is not a complete order: there exist $\l,
\mu$ that can not be compared with respect to this order. 

    The application of this order to our construction is
based on the following key proposition.

    \proclaim{Proposition 5.2} Let $\wt=\tau(\l)w\in \Wt,
w\in W$ and $\ahat=\a+k\delta\in R_\wt$. Write $\wt s_{\ahat}
=\tau(\l')w', w'\in W$. Then:

    \roster \item If $\ahat\in R$ \rom{(}i.e., if
$k=0$\rom{)} then  $\l'=\l$.

    \item If $\ahat\notin R$ \rom{(}i.e., $k>0$\rom{)} then $\l'\prec
\l$. \endroster\endproclaim

    This proposition can be easily proved by direct
calculation.

    Now, let us define the notion of leading term. Let $T$
be an operator in $\C_{q,t}[X]$ of the form 

    $$T=\sum_{\l\in P\v, w\in W} g_{\l, w}(X)\tau(\l)w,
    \tag 5.3$$
where $g_{\l, w}$ are some rational functions in $X$. 

    \proclaim{Definition 5.3} Let $T$ be an operator of the
form \rom{(5.3)}. Assume that it can be written in the
following form:

    $$T=\sum_{w} g_{\l_0, w}(X)\tau(\l_0)w+
    \sum_{\l\prec\l_0, w} g_{\l, w}(X)\tau(\l)w$$
for some $\l_0$ such that at least one of $g_{\l_0, w}\ne
0$. Then we  say that $\sum_{w} g_{\l_0,
w}(X)\tau(\l_0)w$ is the leading term of $T$ and denote it
by $<T>$. \endproclaim

    \remark{Remark} Not every operator has a leading term.
Also, it is not true that the leading term of a product is
the product of leading terms. \endremark

\medskip
{\bf Example 5.4} Consider $Y^\l, \l\in P\v_-$. Then it is easy to see
that $R_{\tau(\l)}$ only contains roots $\ahat=\a+k\delta$
with $k>0$. 
 Thus, if we write 

$$Y^\l=\tau(\l)G^-(\a^{(l)})\ldots G^-(\a^{(1)})$$
and write each $G^-(\a)$ as a sum of two terms (see formula (5.2)) then
-- due to Proposition 5.2 and some simple arguments from the theory of
affine Weyl groups -- the leading term of $Y^\l$ can be obtained if we
replace each $G^-(\a)$ in the expression above by 

$$\frac{t_\a X^\a-t_\a^{-1}}{X^\a-1} ,$$
i.e. if we eliminate the part containing $s_\a$. Thus, 

$$\aligned
<Y^\l>&=\tau (\l)\prod_{\a\in R_{\tau(\l)}}\frac{t_\a
X^\a-t_\a^{-1}}{X^\a-1}\\ 
	&=\left(\prod_{\a\in \tau(\l) R_{\tau(\l)}}\frac{t_\a
X^\a-t_\a^{-1}}{X^\a-1}\right) \tau(\l)\endaligned \tag 5.4$$
(note that this is a product of commuting expressions). 

\medskip

More generally, to find the leading term of $T_\wt$ we have to
separate the affine and non-affine roots in $R_{\wt}$. Let us
introduce the notations 

$$\gathered
R_\wt^0=R_\wt\cap R^+=\{\ahat=\a+k\delta\in R_\wt|k=0\}\\
 R_\wt^{>0}=\{\ahat=\a+k\delta\in R_\wt|k>0\}.\endgathered\tag 5.5$$

\proclaim{Lemma 5.5} For every $\wt\in \Wt$ there exists a reduced
expression $\wt=\pi w, \pi\in \Wt, w\in W$ such that $R_\wt^0=R_w$.
\rom{(}Thus, $R^{>0}_\wt=w^{-1}R_\pi$.\rom{)}\endproclaim

\demo{Idea of proof} Take all reduced expressions $\wt=\pi w$ with
$w\in W$ and choose the one with minimal $l(\pi)$; then use
Lemma~3.2(2).   
\enddemo

\proclaim{Corollary} It is possible to choose a reduced expression for
$\wt$ in such a way that the associated with it sequence of roots
$\a^{(i)}$ looks as follows: 

$$\underbrace{\a^{(l)}, \ldots, \a^{(k)}}_{R_\wt^{>0}}, 
 \underbrace{ \a^{(k-1)}, \ldots, \a^{(1)}}_{R_\wt^0}.$$
\endproclaim

\proclaim{Theorem 5.6} 
\roster \item Let $\wt=\pi w$ be as in Lemma 5.5. Then 

$$<T_\wt>=\left(\prod_{\a\in \wt R_{\wt}^{>0}}\frac{t_\a
X^\a-t_\a^{-1}}{X^\a-1} \right)\pi T_w.\tag 5.6$$

\item Let $\l\in P\v$. Then 

$$<Y^\l> =\sum_w g_w(X)\tau(\l)w $$
for some rational functions $g_w$. 

\endroster
\endproclaim

We could have made more precise statements about the leading term of
$Y^\l$, but this is not necessary for our purposes. 

\demo{Proof} 

(1) Write 

$$T_\wt=\wt G(\a^{(l)})\ldots G(\a^{(k)}) G(\a^{(k-1)})\ldots
G(\a^{(1)})$$
where $\a^{(l)},\ldots, \a^{(k)}\in R_\wt^{>0}, 
\a^{(k-1)},\ldots, \a^{(1)}\in R_\wt^{0} $ (see the Corollary above).
Then, due to Proposition~5.2, we have

$$\aligned
<T_\wt>&=\wt \left(\prod_{\a\in  R_{\wt}^{>0}}\frac{t_\a
X^\a-t_\a^{-1}}{X^\a-1} \right)G(\a^{(k-1)})\ldots
G(\a^{(1)})\\
&= \left(\prod_{\a\in \wt R_{\wt}^{>0}}\frac{t_\a
X^\a-t_\a^{-1}}{X^\a-1} \right) \pi w G(\a^{(k-1)})\ldots
G(\a^{(1)})\\
&= \left(\prod_{\a\in \wt R_{\wt}^{>0}}\frac{t_\a
X^\a-t_\a^{-1}}{X^\a-1} \right)\pi T_w\endaligned
$$

(2) Follows from the expression for $Y^\l$ given in the beginning of
this lecture and Proposition 5.2.
\qed\enddemo

\proclaim{Theorem  5.7} The operators 

$$X^\mu T_\wt, \wt\in \Wt$$
are linearly independent.\endproclaim

\demo{Proof} It suffices to check that their leading terms
are linearly independent. As
for them, if $\wt=\tau(\l) w=\pi w', w, w'\in W$ then 

$$\aligned
<X^\mu T_{\tau(\l) w}>=&
X^\mu \prod_{\a\in \wt R_\wt^{>0}} 
\frac{t_\a X^\a-t_\a^{-1}}
	{X^\a-1} 
\pi T_{w'}\\
=&X^\mu \prod_{\a\in \wt R_\wt^{>0}} 
\frac{t_\a X^\a-t_\a^{-1}}
	{X^\a-1} 
\tau(\l)w {w'}^{-1} T_{w'}.\endaligned$$

    Thus, it is easy to see a relation $\sum_{\l, \mu, w} a_{\l,
\mu,w}X^\mu T_{\tau(\l) w}=0$ is possible only if $\sum_{\mu, w} a_{\l,
\mu,w}X^\mu T_{\tau(\l) w}=0$. Since in the decomposition
$\wt=\tau(\l)w =\pi w'$,  $\pi$ depends only on $\l$, linear
independence of $<X^\mu T_{\tau(\l) w}>$ with fixed $\l$
follows from the fact that
$T_w$ are linearly independent over the field of rational functions of
$X$, which is based on the similar theorem for affine Hecke algebra
$\Hhat^X$. \enddemo

\proclaim{Corollary 5.8} 
\roster\item Elements $X^\mu Y^\l T_w,\l\in P\v, \mu\in P,
w\in W$ are linearly independent in $\H$.
\item $\C_q[X]$ is a faithful representation of $\H$. 
\item Subalgebra $\Hhat ^Y\subset \H$ spanned by $Y^\l,
T_w$ is isomorphic to affine Hecke algebra defined in
Lecture~3, and similarly for $\Hhat^X$. 
\endroster\endproclaim

\proclaim{Theorem 5.9} Let $\pi\in P\v$ be a minuscule coweight. Define
$f_\pi=\sum_{w\in W} Y^{w(\pi)}$, and $L_\pi=\Res f_\pi$, where $\Res $ is
defined by \rom{(4.5)}. Then

$$L_\pi= \sum_{w\in W} w\left( \prod\Sb\a\in R\\(\a, \pi)=1\endSb
\frac{t_\a X^\a -t_\a^{-1}}{X^\a-1}\tau(\pi)\right)\tag 5.7$$
\endproclaim

\demo{Proof} Let us calculate the leading term of $L_\pi$. It follows
from the calculations of leading term for $Y^\l$ (Example 5.4 and
Theorem 5.6(2)) that the leading term has the form $g(X)\tau(\pi^-)$,
where $\pi^-$ is  the antidominant coweight lying in the
orbit of $\pi$: $\pi^-\in P\v_-, \pi^-\in W\pi$. 
 It follows from Theorem~5.6(2) that $\tau(\pi^-)$ can only come from
$Y^{\pi^-}$;  using the calculation of the leading term for
antidominant weight (Example~5.4) we see that the coefficient at
$\tau(\pi^-)$ is equal to 

$$|W_\pi| \prod\Sb\a\in R\\(\a, \pi^-)=1\endSb
\frac{t_\a X^\a -t_\a^{-1}}{X^\a-1}\tau(\pi^-),$$
where $W_\pi$ is the
stabilizer of $\pi$ in $W$.

This gives us the leading term of $L_\pi$. Now, since $\pi$ is
minuscule, it is known that there are no dominant weights $\l$ with
$\l<\pi$, and thus $\l\prec \pi^- \iff \l\in W\pi$. Thus, $L_\pi$ only
contains the terms of the form $g_w(X)\tau(w(\pi))$, which can be
easily calculated, since we know one of them (with $\tau(\pi^-)$) and
$L_\pi$ is $W$-invariant. This gives precisely  formula (5.7).\qed\enddemo

Comparing this formula with the expression for Macdonald's difference
operator $D_\pi$ defined in Lecture~2, we see that they coincide up to
a constant factor. Thus, we see that the operators $L_f, f\in \C[Y]^W$
form a commutative algebra of difference operators which includes the
Macdonald's difference operators -- as was promised in Lecture~2.

\newpage
%%%%%%%%%%%%%%%%%%%%%%%%%%%%%%%%%%%%%%%%%%%%%%%%%%%%%%%%%
\head Lecture 6: Macdonald's polynomials revisited\endhead
%%%%%%%%%%%%%%%%%%%%%%%%%%%%%%%%%%%%%%%%%%%%%%%%%%%%%%%%%%

    As usual, we start with recollections of some results
of previous lectures. We have defined double affine Hecke
algebra $\H$, generated by $X^\mu, \mu\in P, Y^\l, \l\in P\v,
T_w, w\in W$ and defined its action in the space $\C_{q,t}[X]$.
Moreover, we have proved that if $f\in \C_{q,t}[Y]^W\subset \H$  
then the corresponding operator preserves the space
$\C_{q,t}[X]^W$ and its restriction to this space equals 
some $W$-invariant difference operator $L_f$. Also, we
checked that $L_f$ commute and that if $\pi$ is minuscule
coweight, $f=\sum_w Y^{w(\pi)} $ then $L_f$ is Macdonald's
difference operator defined in Lecture~2. In this lecture
we will prove that Macdonald's polynomials are
eigenfunctions of $L_f$ for any $f\in \C_{q,t}[Y]^W$. 
For simplicity, from now on we assume that $t_\a=q^{k_\a}$
for some $k_\a\in \Z_+$; thus, the field $\C_{q,t}$,
considered in the previous lectures becomes $\C_q$.

    Let us start with proving that $L_f$ are triangular in
the basis of $m_\l$.

    \proclaim{Definition} Define a partial order on $P$ as
follows: $\l\prec \mu$ if $\l^+<\mu^+$ or $\l^+=\mu^+$ and
$\l<\mu$, where, as before, $\l^+$ is the dominant weight
lying in the orbit of $\l$. \endproclaim

    Note that this order differs from the order on $P\v$
which we used in the previous lecture; unfortunately, we
have to denote it by the same symbol (there are not so many
symbols available...); we hope it won't cause confusion,
since the order defined in the previous lecture will not be
used in the remaining part of the course. 

\proclaim{Lemma 6.1} Let $\l\in P\v$. Then 

    $$Y^\l X^\mu=\sum_{\nu\preceq \mu} c_{\mu\nu}X^\nu.$$

    In particular, if $\mu\in P^+$ then
$c_{\mu\mu}=q^{2(\l, \mu+\rho_k)}$. \endproclaim

    (Recall that $\rho_k=\frac{1}{2} \sum_{\a\in R^+}
k_\a \a$.)

    \demo{Proof} Assume first that $\l\in P\v_+$. Then the statement
of the Lemma follows from the following two facts, 
 which can be
verified by direct calculation:

    (1) $Y^\l=\tau(\l) G(\a^{(l)})\ldots G(\a^{(1)})$, (see
Lecture~5), where $\a^{(i)}$ run through the set
$R_{\tau(\l)} =\{\ahat=\a+k\delta| \a\in R^+, 0\le k< (\l,
\a)\}$. 

    (2) If $\ahat=\a+k\delta$ is such that $\a\in R^+$ then

    $$G(\ahat)X^\mu=\cases 
t_\a X^\mu +\ldots \quad\text{ if } (\mu, \a\v)\ge 0\\
 t_\a^{-1} X^\mu +\ldots \quad\text{ if } (\mu, \a\v)<0 \endcases
$$

    where dots stand for linear combination of $X^\nu$ with
$\nu\prec \mu$. 

If $\l$ is not dominant, we can write $Y^\l=Y^\mu (Y^\nu)^{-1},
\mu,\nu\in P\v_+$. Since inverse of a triangular matrix is also
triangular, the statement for $\l$ follows from the statements for
$\mu, \nu$. 
\qed
\enddemo

     Lemma 6.1 immediately gives the triangularity of
$L_f$. As before, let $m_\mu=\sum_{\nu\in W\mu} X^\nu$ for
$\mu\in P^+$ be the basis of orbitsums in $\C[X]^W$ 
(we do not consider $m_\mu$ for non-dominant
$\mu$, so whenever a formula contains $m_\mu$ it is always
assumed that $\mu \in P^+$). 

    \proclaim{Theorem 6.2} Let $f\in \C_q[Y]^W, 
\mu\in P^+$. Then 

    $$L_f m_\mu=f(q^{2(\mu+\rho_k)})m_\mu +\sum_{\nu< \mu}
 a_{\mu\nu}m_\nu,\tag 6.1$$
where, by definition, $f(q^\mu)$ is the polynomial in $q$
 obtained by replacing every $Y^\l$ in the
expression for $f$ by $q^{(\l, \mu)}$. 

    \endproclaim

    Thus, it makes sense to talk about the eigenfunctions
of $L_f$. Since a dominant weight $\mu$ is uniquely
determined by the values $f(q^\mu)$ for all $f\in \C[Y]^W$,
it is easy to see that for every dominant $\mu$ there
exists a unique common eigenfunction of $L_f$ in
$\C_q[X]^W$ with the highest term $X^\mu$. Later we will
show that these eigenfunctions are nothing but Macdonald's
polynomials; this gives a new, uniform  proof of the existence of
Macdonald's polynomials. 

    To prove that eigenfunctions are Macdonald's
polynomials, we must check that they are orthogonal with
respect to Macdonald's inner product. Recall that it was
defined in Lecture~2 as follows: for $f, g\in \C_q[X]$ we
let 

    $$\<f,g\>_k=\frac{1}{|W|} [f \bar g \Delta_k]_0,\tag
6.2$$ 
where ($q$-linear) bar involution is defined by
$\overline{X^\mu}=X^{-\mu}$, $[\,\,\,]_0$ is the constant
term, and 

    $$\Delta_k=\prod_{\a\in R} \prod_{i=0}^{k_\a-1}
(1-q^{2i}X^\a).\tag 6.3$$

    This inner product is non-degenerate, $q$-linear,
symmetric and $W$-invariant. 

    However, it turns out that we need to modify this inner
product. Let us introduce the following involution in
$\C_q$: $q^\iota=q^{-1}$ and extend it to $\C_q[X]$,
 letting $(X^\mu)^\iota=X^\mu$. Define
Cherednik's inner product by 

    $$\<f, g\>'_k=[f \bar g^\iota \mu_k]_0,\tag 6.4$$
where

    $$\mu_k= \prod_{\a\in R^+}\prod_{i=1-k_\a}^{k_\a}
(q^iX^{\a/2}-q^{-i} X^{-\a/2}).\tag 6.5$$

    This inner product is not symmetric and not
$W$-invariant. However, it turns out that it is precisely the inner
product suited for our needs, which will become clear very soon. 
Note that the weight function $\mu$ defined above is rather close to
Macdonald's weight function. More precisely, 

$$\mu_k=(-1)^{\sum k_\a} q^{-\sum k_\a(k_\a-1)}\Delta_k 
\frac{\varphi_k}{\delta},\tag 6.6$$
where as before, $\delta=\prod_{\a\in R^+} (X^{\a/2}-X^{-\a/2})$, and 
$\varphi_k=\prod_{\a\in R^+} (q^{k_\a} X^{\a/2}- q^{-k_\a} X^{-\a/2})$.
Note also that $\bar \mu_k=\mu_k^{\iota}$. 

\proclaim{Proposition 6.3} 
\roster\item $\<f, g\>'=(\<g, f\>')^\iota$.

\item If $f, g\in \C_q[X]^W$ then 

$$\<f, g\>'_k= (-1)^{\sum k_\a} q^{-\sum k_\a(k_\a-1)}
  d_k \<f, g^\iota\>_k, \tag 6.7$$
where 

$$d_k =  q^{\sum k_\a}
    \sum_{w\in W} q^{-2\sum_{\a\in R_w} k_\a}.$$
    
\endroster\endproclaim

\demo{Proof} (1) is trivial in view of $\bar \mu_k^\iota=\mu_k$. To
prove (2), note that $[f\bar g^\iota \mu]_0=\frac{1}{|W|} \sum_w[f\bar
g^\iota w(\mu)]_0$, and the result follows from the following
identity: 

    $$\sum_{w\in W} w(\varphi_k/\delta)=\sum_w \prod_{\a\in
R^+}q^{\pm k_\a} =d_k, \tag 6.8$$
where we take the sign $+$ if $w(\a)\in R^+$ and $-$ otherwise. This
identity can be proved in a standard way, by considering the highest
term. \enddemo

\remark{Remark} It is known (see \cite{M4}) that $d_k$ can be written
in the following form:

$$d_k=\prod_{\a\in R^+} 
\frac{q^{(\a\v, \rho_k)+k_\a} -q^{-((\a\v,\rho_k)+k_\a)}}
	{q^{(\a\v, \rho_k)} -q^{-(\a\v,\rho_k)}}.\tag 6.9
$$

However, we are not going to use this formula. 
\endremark

Now we are able to describe 
Macdonald's polynomials in terms of $\<\, ,\, \>'$:

\proclaim{Theorem 6.4} \roster\item $P_\l^\iota=P_\l$.
\item $\<P_\l, P_\mu\>'_k=0$ if $\l \ne \mu$.
\item Macdonald's polynomials are uniquely defined by 
the triangularity condition $P_\l=m_\l+\text{\rm lower
terms}$ and orthogonality 
    condition \rom{(2)} above.
\endroster
\endproclaim

\demo{Proof} (1) Since $\Delta^\iota=const \Delta$, we see that
$[P_\l^\iota \bar P_\mu^\iota \Delta ]_0=0$ for $\mu\ne \l$, and thus
$P_\l^\iota$ satisfy the definition of Macdonald's polynomials. 

(2) follows from (1) and the previous proposition; (3) is obvious since
Cherednik's inner product is non-degenerate. 

\enddemo

Let us now define the notion of adjoint operator. Let $h$ be an
operator in $\C_q[X]$; define its adjoint $h^*$ by the condition 

$$\<hf, g\>'_k=\<f, h^* g\>'_k.$$

An effective way to calculate adjoints is the following. Define a
simpler involution $h\mapsto h^\dag$ by

$$[h(f) \bar g^\iota]_0=[f\, \overline{h^\dag(g)}^\iota]_0;$$ 
thus, $h^\dag$ is the adjoint to $h$ with respect to the inner
product $\<\, ,\, \>'_0$. This adjoint is relatively easy to
calculate; in particular, 

$$\gathered 
p(X)\in \C_q[X]\implies p^\dag=\bar p^\iota\\
\wt \in \Wt \implies \wt^\dag = \wt^{-1}\endgathered$$
(this last condition justifies the introduction of  $\iota$ in the
definition of the inner product: otherwise it would not hold for $\tau(\l)$). 

On the other hand, these involutions are related by a simple rule 
$$h^*=\mu_k ^{-1}h^\dag  \mu_k,$$
 which obviously follows from the definition. In particular, this implies
that $p^*=\bar p^\iota$ for $p\in \C_q[X]$.

\proclaim{Theorem 6.5} \roster 

\item $T_i^*=T_i^{-1}$.

\item $(Y^\l)^*=Y^{-\l}$.
\endroster\endproclaim

\demo{Proof} 

(1) Since $T_i^{-1}=T_i+(t_i^{-1}-t_i), t_i^*=t_i^{-1}$, it suffices
to prove that $(T_i-t_i)^*=T_i-t_i$. From the definition of the action
of $T_i$ we get by direct calculation that 

$$T_i-t_i= \frac{t_iX^{-\a_i/2} -t_i^{-1}X^{\a_i/2}}
		{X^{-\a_i/2}-X^{\a_i/2}}
	(s_i-1).\tag 6.10$$

Since 

$$s_i^*= \mu^{-1} s_i \mu = -\varphi_k^{-1} s_i \varphi_k =
\frac {t_i^{-1}X^{\a_i/2} -t_iX^{-\a_i/2}}
	{t_iX^{\a_i/2} -t_i^{-1}X^{-\a_i/2}}s_i,$$
we get 

$$\aligned
(T_i - t_i)^*=& (s_i^*-1) \frac{t_i^{-1}X^{\a_i/2} -t_iX^{-\a_i/2}}
		{X^{\a_i/2}-X^{-\a_i/2}}\\
	=& \frac{t_i^{-1}X^{\a_i/2} -t_iX^{-\a_i/2}}
		{X^{\a_i/2}-X^{-\a_i/2}} (s_i-1) = T_i- t_i.\endaligned$$ 

(2) It suffices to prove it for $\l\in P\v_+$, in which case it follows
from the previous statement and $\pi_r^*=\pi_r^{-1}$, which can be
proved straightforwardly. 

\qed\enddemo

This means that we can consider ${}^*$ as an involution in $\C_q[Y]$
which is defined by $(f(q)Y^\l)^*=f^\iota Y^{-\l}$; in particular, it
preserves $\C_q[Y]^W$.

\proclaim{Theorem 6.6} Macdonald's polynomials are eigenfunctions of
$L_f, f\in \C_q[Y]^W $: 

$$L_f P_\l=f(q^{2(\l + \rho_k)})P_\l.$$
\endproclaim

\demo{Proof} It is an easy corollary of
$L_f|_{\C_q[X]^W}=f|_{\C_q[X]^W}$ and  Theorems
6.4, 6.5.\enddemo

\proclaim{Proposition 6.7} For $f\in \C_q[Y]^W$,  operators $L_f$ are
self-adjoint with  respect to Macdonald's inner product.\endproclaim

\demo{Proof} It suffices to prove that their restrictions to
$\C_q[X]^W$ are self-adjoint, which follows from the fact that they
are diagonalized in the basis of Macdonald's polynomials, which is
orthogonal with respect to Macdonald's inner product. 
\enddemo

This completes a large part of this course: we have constructed
commuting family of difference operators, whose eigenfunctions are
Macdonald's polynomials. In the next lecture we will apply
this construction to prove the inner product identities. 

\newpage
%%%%%%%%%%%%%%%%%%%%%%%%%%%%%%%%%%%%%%%%%%%%%%%%%%%%%%%%%%%%%%%%
\head Lecture 7: Proof of Macdonald's inner product identities\endhead
%%%%%%%%%%%%%%%%%%%%%%%%%%%%%%%%%%%%%%%%%%%%%%%%%%%%%%%%%%%%%%%%

Recall that we defined action of double affine Hecke algebra $\H$ in
the space $\C_q[X]$. Also, we have defined Cherednik's inner product
$\<\, ,\, \>'_k$ in 
$\C_q[X]$ such that with respect to this inner product
$(Y^\l)^*=Y^{-\l}$, and on symmetric functions it coincides
up to a factor with Macdonald's inner product $\<\, ,\,\>_k$.

The  main goal of this lecture is to prove Macdonald's 
inner product identities (see Theorem 2.4) using the action
of double affine Hecke algebra. From now on, we assume for
simplicity that all $k_a$ are equal: $k_\a=k$, so all $t_\a=t=q^k$.
In fact, it is not much
more difficult to repeat all the arguments for general case; later we
will outline the necessary changes. 

The proof is due to Cherednik; in today's lecture we follow
Macdonald's exposition (\cite{M5}), which simplifies the original
arguments of Cherednik: for example, the introduction of the operator
$\Gt$ below is due to Macdonald.
  
 The idea of proof is quite
simple. First note that due to Proposition~6.3,
  calculation of $\<P_\l, P_\l\>_k$
is equivalent to calculation of $\<P_\l, P_\l\>'_k$. 
Using the large set of operators we have constructed, we want
to prove the theorem by induction in $k$. Let us write $P_\l^{(k)}$ to
denote the dependence of Macdonald's polynomials on $k$. 
We want to construct some
operator $G: \C_q[X]^W\to \C_q[X]^W$ (shift operator), which would
shift $k\to k+1$. More precisely, we want:

 (1) $G P_{\l+\rho}^{(k)}= const P_\l^{(k+1)}$ for some easily computable
constant.

(2) $\<Gf, g\>'_{k+1}=\<f, \Gt g\>'_k$ for some operator $\Gt$. 

In this case, calculation of $\<P_\l^{(k+1)}, P_\l^{(k+1)}\>'_{k+1}$ can be
reduced to calculation of $\<P_{\l+\rho}^{(k)}, \Gt G
P_{\l+\rho}^{(k)}\>'_k$,
 or -- if we know the diagonal terms of $\Gt G$ -- to calculation of
$\<P_{\l+\rho}^{(k)},  P_{\l+\rho}^{(k)}\>'_k$. Repeating the process, we
reduce the question to calculation of $\<P_{\l+k\rho}^{(0)},
P_{\l+k\rho}^{(0)} \>'_0$, which is trivial.

In the $q=1$ case, the shift operators were introduced by Opdam (see
\cite{O3, H3}. The construction for arbitrary $q$ described below is
due to Cherednik (\cite{C6}). 

To define $G,\Gt$, we need the following operators:

$$\gathered
\X = \varphi_{-k} = 
	\prod_{\a\in R^+}(q^{-k} X^{\a/2} -q^{k} X^{-\a/2})\\
\Y= \varphi\v_{-k}(Y) =\prod_{\a\in R^+}(q^{-k} Y^{\a\v/2} -q^{k}
Y^{-\a\v/2})\\
\Yt=\varphi\v_k(Y)= \prod_{\a\in R^+} (q^{k} Y^{\a\v/2}- q^{-k}
Y^{-\a\v/2}).
    \endgathered\tag 7.1$$

It is easily seen from the previous results that 
$\X^\iota=(-1)^{|R^+|}\bar \X=\varphi_k, \X^*=(-1)^{|R^+|}\X$ and 
$\Y^*=(-1)^{|R^+|}\Y, \Yt^*=(-1)^{|R^+|}\Yt$. 

Now, define the shift operators by

$$G=\X^{-1}\Y, \quad \Gt=\Yt \X.\tag 7.2$$

\proclaim{Theorem 7.1} $G, \Gt$ preserve $\C_q[X]^W$. \endproclaim

    \demo{Proof} Recall that $f\in \C_q[X]^W\iff
(T_i-t_i)f=0$ for all $i$ (this follows, for example, from
formula (6.10)). Define 
$\C_q[X]^{-W}=\{f\in \C_q[X]| (T_i+t_i^{-1})f=0\}$. It is easy to see
that as $q\to 1$, this definition becomes the usual definition of
antiinvariant functions. Now, to prove the theorem it suffices to
prove that $ \X(\C_q[X]^W)=
\C_q[X]^{-W}, \Y(\C_q[X]^W)\subset \C_q[X]^{-W},
 \Yt(\C_q[X]^{-W})\subset \C_q[X]^{W}$.
In fact, once we prove that $\X(\C_q[X]^W)\subset
\C_q[X]^{-W}$, the statement that it is isomorphism can be easily
proved by deformation arguments, since in the limit $q\to 1$ this
statement is well-known. Thus, the theorem follows from the following
lemma: 

\proclaim{Lemma}

$$\gathered
(T_i+t^{-1})\X= \frac {t^{-1}  X^{-\a_i/2} -t X^{\a_i/2}}
			{t^{-1}  X^{\a_i/2} -t X^{-\a_i/2}}
		\X (T_i-t)\\
(T_i+t^{-1})\Y= \frac {t^{-1}  Y^{-\a\v_i/2} -t Y^{\a\v_i/2}}
			{t^{-1}  Y^{\a\v_i/2} -t Y^{-\a\v_i/2}}
		\Y (T_i-t)\\
(T_i-t)\Yt = \frac{tY^{-\a_i\v/2} -t^{-1} Y^{\a_i\v/2}}
		{tY^{\a_i\v/2} -t^{-1} Y^{-\a_i\v/2}}	
\Yt (T_i+t^{-1}).
    \endgathered$$
\endproclaim

This lemma is proved by direct calculation, using the identity 
$T_i X^{\a_i/2}-X^{-\a_i/2}T_i = (t-t^{-1})X^{\a_i/2}$. 
\qed
\enddemo

Now we can formulate the main property of shift operators.

\proclaim{Theorem 7.2} Let $f, g\in \C_q[X]^W$. Then 
    
$$\<Gf, g\>'_{k+1} =\frac{d_{k+1}}{d_k} \<f, \Gt g\>'_k,\tag 7.3$$
where $d_k$ are defined in Proposition~6.3.
\endproclaim

\remark{Remark} Since $\mu_{k+1} =\varphi_{k+1}\varphi_{-k} \mu_k$, it is
easy to see that $\<Gf, g\>'_{k+1} =\<f, \Y\varphi_{k+1} g\>'_k$. But
this is of little use, since $\Y\varphi_{k+1}$ does not preserve
$\C_q[X]^W$; thus, to find, say, $\<P_\l, \Y\varphi_{k+1} m_\mu\>'_k$ we
have to calculate the highest term of projection of $\Y\varphi_{k+1}
m_\mu$ on $\C_q[X]^W$, which is very difficult. \endremark

\demo{Proof} The proof is based on the following simple idea, which
we have already used before. Let $\P =\frac{1}{|W|} \sum_{w\in W} w$ be
the usual symmetrizer. Then for every $f\in \C_q[X]$ we have
$[f]_0=[\P f]_0$. Thus, if $\P f= \P g$ then $[f]_0=[g]_0$. Also, we
will need the following proposition. 

\proclaim{Proposition 7.3} Let
 $\P_- =\frac{1}{|W|} \sum_{w\in W}(-1)^{l(w)} w$
be the usual antisymmetrizer. Then for every $f\in \C_q[X]^W$ we have 

$$\P_- \Y f = \P_- \Yt f.$$
\endproclaim

The proof of this proposition is quite non-trivial and requires
introduction of new interesting operator -- $q$-antisymmetrizer. We
will give this proof in the next lecture. 

Now let us prove the theorem. By definition,

$$
\<Gf, g\>'_{k+1} =[(Gf) \bar g^{\iota}\mu_{k+1}]_0
 =[(Gf) \bar g^{\iota}\P(\mu_{k+1})]. $$ 

    Since $\mu_{k+1}=\varphi_{k+1}\varphi_{-k} \mu_k=\varphi_{k+1}\X
\mu_k$, and $\X\mu_k$ is antisymmetric (this is the crucial
step!), we have $\P(\mu_{k+1})=
\P_-(\varphi_{k+1})\X\mu_k=\frac{1}{|W|} d_{k+1} \delta \X\mu_k$ (see
formula (6.8)). 

    Similarly, $\P (\X^2\mu_k)=\frac{1}{|W|}d_k\delta \X\mu_k$, and
thus, 
$$\P(\mu_{k+1})=\frac{d_{k+1}}{d_k} \P(\X^2 \mu_k).$$

    Substituting it in the expression for $\<Gf,
g\>'_{k+1}$, we get 

    $$\split
\<Gf,& g\>'_{k+1}=\frac{d_{k+1}}{d_k} [\X^{-1}\Y (f)
\bar g^{\iota} \X^2\mu_k]_0\\
&=\frac{d_{k+1}}{d_k} [\P(\Y (f)
\bar g^{\iota} \X\mu_k)]_0
=\frac{d_{k+1}}{d_k} [\P_-(\Y (f))
\bar g^{\iota} \X\mu_k]_0.\endsplit$$

Using Proposition 7.3 we can replace in the last formula $\Y$ by 
$\Yt$, and thus

    $$\split\<Gf,& g\>'_{k+1} =\frac{d_{k+1}}{d_k} [\X\Yt (f)
\bar g^{\iota} \mu_k]_0\\
&=\frac{d_{k+1}}{d_k} \<\X\Yt f, g\>'_k
=\frac{d_{k+1}}{d_k} \<f, \Yt\X g\>'_k.\endsplit$$
\qed
\enddemo

    This immediately implies that the shift operators
indeed shift the parameter $k$ of Macdonald's polynomials:

    \proclaim{Theorem 7.4} 
\roster\item $G P_{\l+\rho}^{(k)}=q^{k|R^+|}c_k (\l) P_\l^{(k+1)}$,
where 

    $$c_k(\l)=\prod_{\a\in R^+} 
(q^{-k+(\a\v, \l+(k+1)\rho)}-q^{k-(\a\v, \l+(k+1)\rho)}),\tag
7.4$$

    and $G P_\l^{(k)}=0$ if $\l-\rho\notin P_+$.
\item For $k\ge 0$, $\Gt P_\l^{(k+1)} = 
q^{-k|R^+|}\hat c_k(\l) P_{\l+\rho}^{(k)}$, where 

    $$\hat c_k(\l)=\prod_{\a\in R^+} 
(q^{k+(\a\v, \l+(k+1)\rho)}-q^{-k-(\a\v, \l+(k+1)\rho)}).\tag
7.5$$
\endroster\endproclaim

    \demo{Proof} First, it is easy to prove, using Lemma~6.1
that $G P_{\l+\rho}^{(k)}=q^{k|R^+|}c_k m_\l+\ldots$. Thus,
to prove (1) it suffices to check that
$\<GP_{\l+\rho}^{(k)}, m_\mu\>'_{k+1}=0$ if $\mu<\l$. Due to
Theorem~7.2, this is equivalent to $\<P_{\l+\rho}^{(k)},
\Gt m_\mu\>'_k=0$. Since $\Gt m_\mu$ is linear combination
of $m_\nu$ with $\nu\le \mu+\rho$ (this also follows from
Lemma~6.1), the statement follows from the definition of
Macdonald's polynomials. (2) is proved in a similar
way.\qed\enddemo

Now we can prove Macdonald's inner product identities. 
Let us introduce $M'_k(\l)=\<P_\l^{(k)}, P_\l^{(k)}\>'_k$. 

\proclaim{Proposition 7.5} 

    $$M'_{k+1}(\l)=(-1)^{|R^+|}\frac{d_{k+1}}{d_k}
		\frac{\hat c_k(\l)}{c_k(\l)} M'_k(\l+\rho).\tag 7.6$$
  \endproclaim

 \demo{Proof} Using the previous theorem, we can write

    $$\split
M'_{k+1}(\l)=& (c_k(\l) c_k(\l)^\iota)^{-1}
\<GP_{\l+\rho}^{(k)}, GP_{\l+\rho}^{(k)} \>'_{k+1}\\ 
= &\frac{d_{k+1}}{d_k}(c_k(\l) c_k(\l)^\iota)^{-1}
	\<P_{\l+\rho}^{(k)}, \Gt G P_{\l+\rho}^{(k)}\>'_k\\
=& \frac{(\hat c_k(\l)c_k(\l))^\iota}{c_k(\l) c_k(\l)^\iota} 
\frac{d_{k+1}}{d_k}\<P_{\l+\rho}^{(k)}, P_{\l+\rho}^{(k)}\>'_k    
=(-1)^{|R^+|}\frac{d_{k+1}}{d_k}\frac{\hat c_k(\l)}{c_k(\l)}
M'_k(\l+\rho).
\endsplit$$
\qed\enddemo
    
    \proclaim{Corollary 7.6} Let $M_k(\l)=\<P_\l^{(k)},
P_\l^{(k)}\>_k= d_k^{-1} (-1)^{k|R^+|} q^{k(k-1)|R^+|}
M'_k(\l)$ \rom{(}see Proposition~\rom{6.3)}. Then 
$$M_{k+1}(\l)=\prod_{\a\in R^+} 
\frac{1-q^{2(\a\v, \l+(k+1)\rho)+2k}}
{1-q^{2(\a\v, \l+(k+1)\rho)-2k}} M_k(\l+\rho).\tag 7.7$$ 
\endproclaim

    Applying this corollary $k-1$ times and using
$M_1(\l)=1$ for all $\l$ (this is equivalent to saying that
Weyl characters are orthonormal), we get Macdonald's inner
product identities, formulated in Lecture~2.
 \proclaim{Theorem 7.7}{\rm (Macdonald's inner product
identities)} If all $k_\a=k$ then  
$$
\<P_\l, P_\l\>_k=\prod_{\a\in R^+}
\prod_{i=1}^{k-1} 
\frac{1-q^{2(\a\v, \l+k\rho)+2i}}
    {1-q^{2(\a\v, \l+k\rho)-2i}}.$$
 \endproclaim

    To prove inner product identities in general case, i.e.
when $k_\a$ are not necessarily equal (see Theorem~2.4), we
have to introduce shift operators separately for
long and short roots. They are defined in precisely the same way as we
did, but with product in (7.1) only over long (respectively, short)
roots. Repeating the steps above with necessary changes, we can prove
that these shift operators change $k_\a$ for long (respectively,
short) roots by one, and prove general Macdonald's inner product
identities (2.6). We refer the reader to \cite{C6} for details.

\newpage
%%%%%%%%%%%%%%%%%%%%%%%%%%%%%%%%%%%%%%%%%%%%%%%%%%%%%%%%%%%%%%
\head Lecture 8: $q$-symmetrizers \endhead
%%%%%%%%%%%%%%%%%%%%%%%%%%%%%%%%%%%%%%%%%%%%%%%%%%%%%%%%%%%%%%%%%%%%

In this lecture we prove 
 Proposition~7.3 and thus complete the proof of inner product
identities. Recall that we want
to prove $\P_-\Y f=\P_- \Yt f$ for every $f\in \C_q[X]^w$, where
$\P_-$ is the antisymmetrizer. 
Unfortunately, commutation relations of $w\in W$ with $Y$ are very
complicated, which makes direct calculation impossible. However, there
is a bypass, which involves introduction of $q$-antisymmetrizer; this
does not seem to be closely related with Macdonald's theory, but is
interesting enough in itself, so we spend some time discussing these new
operators.

Let us start with describing of kernel of the antisymmetrizer. 

\proclaim{Theorem 8.1} In any  finite-dimensional representation $V$ of $W$
we have 

$$\Ker \P_- =\sum_i \Ker (1-s_i).\tag 8.1$$
\endproclaim

\demo{Proof} It is clear that $\Ker(1-s_i) \subset \Ker \P_-$, so the
difficult part is to prove equality. If $V$ is a representation of
$W$, denote $V_i=\Ker (1-s_i), V'=\sum V_i$. 

\proclaim{Lemma} $V'$ is $W$-invariant.
\endproclaim

\demo{Proof} It suffices to prove $s_i V_j\subset V_i+V_j$. Let $v\in
s_i V_j$; then $s_j (s_i v)=s_i v$. Introduce $v_{\pm}
=\frac{1}{2}(v\pm s_iv)$. Then $v=v_+ +v_-, s_i v=v_+ - v_-$, and thus
$s_j (v_+ - v_-)= v_+ - v_-$, so $v_+ - v_-\in V_j$. Since by
definition $v_+\in V_i$, we see that $v\in V_i + V_j$.\enddemo

Now, let us return to the proof of the theorem. Obviously, it suffices
to prove this theorem for an irreducible representation. In this case,
$V'$ can be either $0$ or $V$. But: 

$$ \gathered
V'=0\iff
\text {all } V_i=0\iff\\
\text{ for all }i, \quad (1-s_i) \text{ is invertible} \iff\\
\text{ for all }i, \quad s_i=-1 \text{ in } V \iff\\
V \text{ is the sign representation.}\endgathered$$

Thus, for an irreducible $V$ we have

$$V'=\cases 0,\quad  V \text{ is the sign representation}\\
	    V\quad \text{ otherwise}\endcases$$

Obviously, this coincides with $\Ker \P_-$. \qed\enddemo

Note that (8.1) also holds for  the representation of $W$ in the space
of polynomials $\C_q[X]$, since this representation is a direct sum of
finite-dimensional representations. 

The main idea of proof of Proposition~7.3 
is that now  we can describe
$\Ker \P_-$ in $\C_q[X]$ 
using the action of the Hecke algebra $H$  generated by
$T_1, \dots, T_n$ rather then the action of $W$, and  then use the
commutation relations of $H$ with $Y$.

Let us introduce the following element of $H$ which we will call the 
$q$-antisymmetrizer: 

$$\P^q_-= d^{-1} \sum_{w\in W} (-t)^{-l(w)} T_w,\tag 8.2$$
    where $d=\sum_{w\in W} t^{-2 l(w)}$. 

It is easy to see that as $q\to 1$ this element becomes the usual
antisymmetrizer $\P_-$. 

\proclaim{Theorem 8.2} 
\roster\item For every $i=1,\dots, n$, $\P^q_-$ is divisible by $T_i-t$ both on
the left and on the right. 

\item We have the following properties for the action of $\P^q_-$ in
$\C_q[X]$ :

$$\gathered\Ker  \P^q_- =\Ker \P_-\\
\operatorname{Im } \P^q_- = \C_q[X]^{-W}.\endgathered$$

\item $\P^q_-$ is a projector. 
\endroster
\endproclaim

\demo{Proof} (1) Since $w\mapsto ws_i$ is an involution of the Weyl
group, $W$ is a
union of pairs $w, ws_i$ where $w$ is such that $l(ws_i)=l(w)+1$.
Thus,

$$\P^q_-=d^{-1}\sum_{l(ws_i)=l(w)+1}(-t)^{-l(w)}T_w(1-t^{-1}T_i).$$

Divisibility on the left is proved similarly.

(2) It follows from (1)  that $\Ker \P^q_- \supset \sum
\Ker (T_i-t)=\sum \Ker (s_i-1)$. On the
other hand Theorem~8.1 claims  that for $q=1$ this inclusion is an
equality. Since the rank can not increase under specialization, it
implies  $\Ker \P^q_- = \Ker \P_-=\sum \Ker (1-s_i)$. 
Similarly, (1) implies that
 $\operatorname{Im} \P^q_- \subset
\C_q[X]^{-W}$; since the dimension of $\C_q[X]^{-W}$ is the same as
for $q=1$, we see that it is also an equality. (Of course, to make
sense of these dimension arguments we must consider $\C_q[X]$ as a
filtered space and note that both $\P_-, \P^q_-$ preserve this filtration.)

    (3) is quite trivial: let $v\in \C_q[X]^{-W}$. Then
$T_iv =-t^{-1} v$, so $T_w v = (-t)^{-l(w)} v$, and $\P_-^q
v =v$. \qed
\enddemo

\proclaim{Corollary} For $f\in \C_q[X], \P_- f=0\iff \P_-^q
f=0$.\endproclaim

Thus, to prove Proposition~7.3 it suffices to prove 
$\P^q_-(\Y-\Yt)f =0$ for every $f\in \C_q[X]^W$.
 Using the fact that $C_q[X]$ is a faithful
representation of $\H$, it is easy to prove that this last condition
is equivalent to 

$$ \P^q_-(\Y-\Yt)= \sum h_i (T_i-t) \quad\text{for some } h_i\in
\Hhat^Y\tag 8.3$$
as elements of $\Hhat^Y$. 

Now, we can do the following trick. Since (8.3) is an identity in
$\Hhat^Y$, it suffices to prove it in any faithful representation of
$\Hhat^Y$. Let us prove it in $\C_q[Y]$ (see Lecture~3). Now we can
repeat the same chain of arguments in the reverse order: (8.3) $\iff
\P_-^q(\Y-\Yt) f=0$ for every $f\in \C_q[Y]^W \iff \P_- (\Y-\Yt) f=0$.
But this last condition is trivial: since $\Y, \Yt$ act in $\C_q[Y]$
just by multiplication, action of $W$ is just by permuting indices of
$Y^\l$; in particular,  $w_0(\Y)=(-1)^{|R^+|}\Yt, 
w_0(\Yt)=(-1)^{|R^+|}\Y $, where $w_o$ is the longest
element of the Weyl group. Since $\P_-$ is divisible by 
$(1+(-1)^{|R^+|}w_0)$, it implies $\P_- (\Y-\Yt)=0$. This
completes the proof of Proposition~7.3, and thus, of
Macdonald's inner product identities.


\Refs
\widestnumber\key{AAA}

\ref\key AI \by Askey, R. and Ismail, M.E.H. 
\paper A generalization of ultraspherical polynomials
\inbook Studies in Pure Mathematics
\ed P. Erd\"os\publ Birkh\"auser \yr 1982\pages 55--78\endref

\ref \key AW \by Askey, R. and Wilson, J.
\paper Some basic hypergeometric orthogonal polynomials
that generalize Jacobi polynomials
\jour Memoirs of AMS \vol 319\yr 1985\endref


\ref\key B \by Bourbaki, N.
\book Groupes et alg\`ebres de Lie, Ch. {\bf 4--6}
\publ Hermann \publaddr Paris\yr 1969\endref


\ref \key BZ\by  Bressoud, D. and Zeilberger, D.
\paper  A proof of Andrews' $q$-Dyson conjecture
\jour Discrete Math.\vol 54\yr 1985\pages 201--224\endref

\ref\key C1\by Cherednik, I. 
\paper Double affine Hecke algebras, 
Knizhnik- Za\-mo\-lod\-chi\-kov equa\-tions, and Mac\-do\-nald's 
ope\-ra\-tors 
\jour IMRN (Duke M.J.) \vol   9\yr 1992 \pages 171--180\endref


\ref\key C2\bysame
\paper  The Macdonald constant term conjecture
\jour IMRN  \vol 6 \yr 1993\pages 165--177
\endref 

\ref \key C3\bysame
\paper   A unification of Knizhnik--Zamolodchikov
and Dunkl operators via affine Hecke algebras
\jour Inventiones Math.\vol  106\issue 2\yr 1991\pages 411--432\endref


\ref\key C4\bysame
\paper Quantum Knizhnik--Za\-mo\-lod\-chi\-kov
equa\-tions and affine
root systems \jour Commun. Math. Phys. \vol 150\yr 1992\pages
109--136
\endref

\ref\key C5\bysame
\paper Integration of Quantum many-body problems by affine 
Knizhnik--Za\-mo\-lod\-chi\-kov equations 
\jour Pre\-print RIMS--776 \yr 1991\finalinfo
(Advances in Math.(1994))\endref

\ref\key C6\bysame
\paper Double affine Hecke algebras and Macdonald conjectures
\jour to appear in Annals of Math. \yr 1994\endref

\ref\key C7\bysame
\paper Difference-elliptic operators and root systems\jour 
preprint, November 1994\endref

\ref\key D\by van Diejen, J.F.
\paper Commuting difference operators with polynomial eigenfunctions
\jour to appear in Compos. Math. \endref


\ref\key EFK\by Etingof, P.I., Frenkel, I.B. and Kirillov, A.A.,
Jr\paper Spherical functions on affine Lie groups\jour hep-th 9403168 
(submitted to Duke Math. J)\yr 1994\endref 



\ref \key GG\by Garvan, F. and  Gonnet, G.
\paper Macdonald's constant term conjectures
for exceptional root systems
\jour Bull. AMS\vol 24\issue 2\yr 1991\pages 343--347\endref

\ref\key Ha \by Habsieger, L.
\paper La $q$-conjecture de Macdonald-Morris pour $G_2$
\jour C.R.Acad. Sci. Paris S\'er.{ 1} Math.\vol 303\yr 1986\pages
211-213\endref

\ref \key HO \by Heckman, G.J., Opdam, E.M.\paper Root systems and
hypergeometric functions I\jour Compos. Math.\vol 64 \pages 329--352\yr
1987\endref

\ref \key H1 \by Heckman, G.J.\paper Root systems and
hypergeometric functions II\jour Compos. Math.\vol 64 \pages 353--373\yr
1987\endref

\ref\key H2\bysame
\paper A remark on the Dunkl differential-difference operators
\inbook Harmonic analysis on reductive groups
\eds W. Barker, P. Sally
\publ Birkh\"auser\yr 1991\pages 181--191\endref

\ref \key H3
\bysame
\paper  An elementary approach to the hypergeometric shift operators of
Opdam \jour Invent.Math. \vol  103\yr 1991\pages 341--350\endref


\ref\key Hu1\by Humphreys, J.E. \book Introduction to Lie algebras and
representation theory\publ Springer-Verlag\publaddr New York\yr
1972\endref

\ref\key Hu2 \bysame
\book Reflection groups and Coxeter groups \publ Cambridge Univ.
Press\publaddr Cambridge\yr 1990\endref


\ref\key K\by Kadel, K.
\paper A proof of the q-Macdonald-Morris conjecture
for $BC_n$\jour preprint\endref

\ref\key Ko\by Koornwinder, T.H.
\paper Special functions associated with root systems: recent progress
\inbook From Universal Morphisms to Megabytes --- a Baayen Space
Odyssey \eds K. R. Apt, A. Schrijver, \& N. M. Temme
\publ CWI, Amsterdam \yr 1994\pages 391--404\endref





\ref\key L
\by  Lusztig, G. 
\paper Affine Hecke algebras and their graded version
\jour J. of the AMS \vol 2\issue 3 \yr 1989\pages 599--685 \endref

\ref\key M1\by Macdonald, I.G. \paper A new class of symmetric
functions\jour Publ. I.R.M.A. Strasbourg, 372/S-20, Actes 20
S\'eminaire Lotharingien\pages 131-171\yr 1988\endref


\ref\key M2\bysame\paper Orthogonal polynomials associated
with root systems\jour preprint\yr 1988\endref

\ref\key M3\bysame
\paper Some conjectures for root systems
\jour SIAM J. of Math. Analysis\vol 13\issue 6\yr 1982
\pages 988--1007\endref

    
\ref\key M4\bysame
\paper The Poincar\'e series of a Coxeter group
\jour Math. Annalen\vol 199\yr 1972\pages 161--174\endref

\ref\key M5\bysame 
\paper Orthogonal polynomials and constant term conjectures
\jour Lectures at Leiden University\yr May 1994\endref 

\ref\key Ma \by Matsuo, A. \paper Integrable connections related to
zonal spherical functions\jour Inv. Math.\vol 110 \pages
95--121 \yr 1992\endref

\ref\key N\by Noumi, M. 
\paper  Macdonald's symmetric polynomials as zonal
spherical functions on quantum homogeneous spaces
\jour Adv. in Math.\yr 1995 \toappear\endref

\ref \key O1 \by Opdam, E.M.\paper Root systems and
hypergeometric functions III\jour Compos. Math.\vol 67\pages 21--49\yr
1988\endref

\ref \key O2 \bysame\paper Root systems and
hypergeometric functions IV\jour Compos. Math.\vol 67 \pages 191--207\yr
1988\endref


\ref\key O3 \bysame
\paper Some applications of hypergeometric shift operators
\jour Inv. Math.\vol 98\yr 1989\pages 1--18\endref

\ref\key OOS\by Ochiai, H., Oshima, T., and Sekiguchi, H.
\paper Commuting families of symmetric differential operators
\jour Proc. of the Japan Acad.\vol 70, Ser. A \issue 2\yr 1994\pages
62--68\endref


\ref\key OP\by Olshanetsky, M.A.  and Perelomov, A.M.\paper Quantum
integrable systems related to Lie algebras \jour Phys. Rep. \vol 94
\pages 313-404\yr 1983
\endref

\ref\key Su \by Sutherland, B.\paper Exact results for quantum many-body
problem in one dimension\jour Phys. Rep. \vol A5 \pages 1375--1376\yr
1972\endref 


\ref\key V\by Verma, D-N.
\paper The role of affine Weyl groups in the 
representation theory of algebraic Chevalley groups
and their Lie algebras
\inbook Lie groups and their representations (Proceedings of 
the Summer School on Group Representations)
\publaddr Budapest \yr 1971  \pages 653--705\endref


\endRefs
    \enddocument
\end